\documentclass[11pt,reqno,a4paper]{amsart}

\usepackage[utf8]{inputenc}
\usepackage{amsfonts,latexsym,amssymb,amsmath,amsthm,amsrefs,tikz}
\usepackage{hyperref}
\hypersetup{pdftitle={Microlocal properties of Cauchy evolution under metric variations},pdfauthor={Alexander Strohmaier}}
\usepackage{a4wide}
\newcommand{\Id}{\mathrm{Id}}
\DeclareMathOperator{\supp}{supp}
\DeclareMathOperator{\WF}{WF}

\newcommand{\der}{\mathrm{d}}
\newcommand{\rmi}{\mathrm{i}}
\newcommand{\R}{\mathbb R}
\newcommand{\C}{\mathbb C}
\newcommand{\calM}{\mathcal M}
\newcommand{\calH}{\mathcal H}
\newcommand{\calD}{\mathcal D}
\newcommand{\calL}{\mathcal L}

\newcommand{\dotT}{\dot T^*}
\newcommand{\calF}{\mathcal F}
\DeclareMathOperator{\Hom}{Hom}
\newtheorem{theorem}{Theorem}

\newtheorem{lemma}[theorem]{Lemma}
\newtheorem{proposition}[theorem]{Proposition}
\newtheorem{corollary}[theorem]{Corollary}
\newtheorem{assumption}{Assumption}

\theoremstyle{remark}
\newtheorem{remark}[theorem]{Remark}

\title[Cauchy evolution under metric variations]{Microlocal properties of Cauchy evolution under metric variations}
\author[A.~Strohmaier]{Alexander Strohmaier}
\address{Leibniz University Hannover, Institute of Analysis, 30167 Hannover, Germany}
\email{a.strohmaier@math.uni-hannover.de}
\date{}

\begin{document}
\begin{abstract}
We develop the microlocal analysis of the Cauchy evolution for wave-type equations on the infinite-dimensional space of Lorentzian metrics. We regard the Cauchy evolution operator as a distributional function of the metric, and give this notion a precise meaning through restrictions to finite-dimensional families of metrics. We study the singularities of this metric dependence.
For normally hyperbolic operators on vector bundles, we determine the resulting operator-valued wavefront sets. Their parameter covectors are given by integrals of metric variations along individual null geodesics.
At ultrastatic ends, metric dependence also determines the time-translation generators and their positive- and negative-energy splitting. More generally, the parameter singularities distinguish the two characteristic sheets and can be used to reconstruct the ordinary Cauchy propagation relation, which describes propagation of singularities.
Consequences include smoothing under parameter averaging and polar-cone bounds associated with averaged null-positive variations, linking the result to averaged null energy conditions. The analysis also covers the Cauchy evolution of the Dirac equation.
\end{abstract}
\maketitle

\section{Introduction and main results}\label{sec:setup}

Let $(M,g_0)$ be a time-oriented globally hyperbolic spacetime with signature $(+1,-1,\ldots,-1)$. Assume $E \to M$ is a complex vector bundle of finite positive rank. A second-order differential operator $P: C^\infty(M;E) \to C^\infty(M;E)$ is called normally hyperbolic if, with our sign convention, its principal symbol is $\sigma_2(P)(x,\xi)=-g_0^{-1}(\xi,\xi)\Id_{E_x}$. Equivalently, the leading second-order term is the same as that of the d'Alembert operator and therefore such operators are also referred to as d'Alembert-type or wave-type operators. Examples are the Klein--Gordon operator and the square of any Dirac-type operator.

For distributional solutions of $Pu=0$, the wavefront set lies in the null characteristic set and is invariant under the null bicharacteristic flow (see for example \cites{Ho1, MR0388464, MR2512677, HormanderPropagation1971}). Since the light cone decomposes into future and past cones, this leads to natural microlocal decompositions of the associated propagators, giving rise to Hadamard states and Feynman propagators.
Hadamard states are of particular importance in quantum field theory on curved spacetimes and have been developed in the physics literature
(see for example \cites{FullingSweenyWaldHadamard,FullingNarcowichWaldHadamard}). The relevance of the microlocal splittings in the earlier theory of
distinguished Feynman parametrices of Duistermaat and H\"ormander
\cite{MR0388464} in this context has been recognised by Radzikowski \cite{RadzikowskiHadamard}, leading to the use of microlocal analysis and wavefront sets in quantum field theory, including the analogous microlocal characterisation for vector-valued fields (\cite{SahlmannVerchHadamard}). 
Hadamard states are non-unique and their construction, retaining positivity properties required for quantum field theory, has inspired several other developments in microlocal analysis. For constructions of Hadamard states and related frequency splittings by pseudodifferential methods see
\cites{JunkerHadamard,JunkerHadamardErratum,JunkerSchroheAdiabatic,GerardWrochnaHadamard};
for constructions from null boundary data see \cite{DappiaggiMorettiPinamontiHadamard},
and for operator-theoretic and global microlocal constructions of Feynman propagators see
\cite{DerezinskiSiemssenFeynman} and \cite{VasyWrochnaGlobal}, respectively. For the relation of bundle-valued
Feynman propagators to Hadamard states see
\cites{IslamStrohmaierFeynman,FewsterStrohmaierHadamard}. Microlocal splittings also underlie the Fredholm and index theory of
Lorentzian Dirac operators with Atiyah--Patodi--Singer boundary
conditions. They imply compactness of the off-diagonal blocks of
Cauchy evolution relative to the boundary spectral projections,
which gives Fredholmness of the APS problem
\cite{BaerStrohmaierIndex}*{Sections 2 and 3}. This connection has been
used in several directions in Lorentzian index theory
\cites{BaerStrohmaierLocalIndex,ShenWrochnaIndex,BaerHannesBoundary,
BravermanCalliasIndex,VanDenDungenRongeAPS}, and more recently in work
on $L^2$-$\Gamma$ indices, equivariant indices, and Feynman conditions
\cites{DamaschkeVertmanIndex,IslamRongeEquivariantIndex,AmarFeynmanIndex}.

The aim of this paper is to show that instead of wavefront sets on the base manifold one can also use metric deformations to characterise these microlocal splittings and  the propagation relation. We do this by studying
the singular dependence of Cauchy evolution on the metric. Regarding metric parameters as microlocal variables, we determine the operator-valued parameter wavefront sets, separately for the two characteristic sheets, and show that the ordinary Cauchy propagation relation can be recovered from them. This gives a formulation of characteristic propagation through metric dependence. At ultrastatic ends, diffeomorphism-induced metric variations also recover the time-translation generators and the related positive- and negative-energy splittings. 
Under a mild non-degeneracy hypothesis (Assumption \ref{ass:A}) integration with a test function smoothens the evolution; positivity of the ray integrals gives associated polar-cone bounds. This characterises the parameter wavefront sets of the microlocal positive- and negative-energy components in terms of metric variations. The polar-cone bounds are analogous to the microlocal spectrum condition for metric variations proposed in \cite{StrohmaierStress}*{Definition 7}. For the relation of metric dependence to
relative Cauchy evolution and the stress-energy tensor see
\cite{MR2007173}, and for analytical constructions of metric-dependent
wave propagation see \cite{StrohmaierStress}.

We will now describe the precise setting and main results.
To keep the exposition simple we assume throughout that $(M,g_0)$ is spatially compact, i.e. its Cauchy surfaces are compact. Note that, for any two Cauchy surfaces in a globally hyperbolic spacetime, the Cauchy evolution operator acting on initial data supported in a fixed compact subset can be described by the Cauchy evolution operator of a spatially compact spacetime. This is a consequence of finite speed of propagation of solutions for normally hyperbolic differential operators. There is therefore a straightforward modification of the main statements of this paper for general globally hyperbolic spacetimes, by letting Cauchy evolution act on the locally convex spaces of compactly supported Cauchy data. We do not carry out this extension here. 

Let $K\subset M$ be compact. We denote by $\calM_K$ the set of smooth globally
hyperbolic Lorentzian metrics $g$ with $\supp(g-g_0)\subset K$.
Let $B$ be an open neighbourhood of $0$ in $\R^N$ and let
\[
 B\longrightarrow\calM_K,\qquad b\longmapsto g_b
\]
be a smooth family in the sense that the map $(b,x)\mapsto g_b(x)$ from $B \times M$ to $\operatorname{Sym}^2T^*M$ is smooth.
For $v\in T_bB$ we write
\[
 h_b(v):=\der_b g(v)\in C^\infty_K(M;\operatorname{Sym}^2T^*M).
\]

Let $P_b:C^\infty(M;E)\to C^\infty(M;E)$ be a smooth family of normally hyperbolic operators. In local coordinates and a bundle trivialisation they have the form
\[
 P_b=g_b^{\mu\nu}\mathbf{1} \partial_\mu\partial_\nu
          +v_b^\mu\partial_\mu+w_b,
\]
where $v_b^\mu$ and $w_b$ are matrix-valued smooth coefficient functions,  $\mathbf{1}$ denotes the identity matrix, and the Einstein summation convention is used.The principal symbol is thus 
\[
 \sigma_2(P_b)(x,\xi)=-g_b^{-1}(\xi,\xi)\Id_{E_x}.
\]
The operator determines a unique connection $\nabla^b$ on $E$ and a unique smooth endomorphism $W_b$ such that, with this sign convention,
\begin{equation}\label{eq:Weitzenboeck}
 P_b=\operatorname{tr}_{g_b}(\nabla^b)^2+W_b.
\end{equation}
We call this the
\emph{Weitzenb\"ock formula} and $\nabla^b$ the \emph{Weitzenb\"ock connection}.
Smoothness of the family $P_b$ means that its coefficients in local bundle trivialisations depend smoothly on $(b,x)$.
The Weitzenb\"ock connection therefore depends smoothly on $b$. All arguments are local in parameter space, so a smoothly varying family of bundles may be identified with a fixed bundle after shrinking $B$. We therefore assume without loss of generality that $E$ is independent of $b$.
We assume that all coefficients of $P_b$ agree
with those of $P_0$ outside $K$. Lower-order terms may otherwise depend
smoothly on $b$. Examples are the scalar operator $P_b=\Box_{g_b}+m^2$ and the square $P_b=D_b^2$ of a Dirac-type operator, with the sign convention for which its principal symbol is the one displayed above.
 The
Hamiltonian used below is
\begin{equation}\label{eq:H}
 H_b(x,\xi)=\frac12g_b^{-1}(\xi,\xi).
\end{equation}
Its null bicharacteristics are independent of the lower-order terms and correspond to the null geodesics of the Lorentzian metric $g_b$. Their projections satisfy $\dot\gamma=\xi^\sharp$, so increasing Hamiltonian parameter propagates future-directed covectors towards the future and past-directed covectors towards the past.

After shrinking $B$ we choose a common Cauchy temporal function $t$ and a splitting $M = \R_t \times \Sigma$. Choose fixed $t_-<t_+$ and let $\Sigma_\pm=\{t_\pm\}\times\Sigma$ be the boundary surfaces of the slab. We do not require $K$ to be disjoint from these surfaces; the coefficients of $P_b$ and the normal vectors may depend on $b$ there. Fix an auxiliary smooth connection $\nabla$ on $E$, independently of $b$, for recording the Cauchy data. It need not be the Weitzenb\"ock connection of $P_0$. For a
solution of $P_bu=0$ its Cauchy data are
\[
 R_{b,\Sigma}u=(u|_\Sigma,\nabla_{n_b}u|_\Sigma),
\]
where $n_b$ is the future unit normal. Write $E_\Sigma=E|_\Sigma$ and $E_\pm=E|_{\Sigma_\pm}$. The Cauchy evolution is
\begin{equation}\label{eq:U}
 U_b=R_{b,\Sigma_+}R_{b,\Sigma_-}^{-1}.
\end{equation}
For $r\in\R$ we use the scale
\[
 \calH^r(\Sigma)=H^{r+\frac12}(\Sigma;E_\Sigma)\oplus H^{r-\frac12}(\Sigma;E_\Sigma).
\]
The Sobolev norms, auxiliary positive definite Hermitian metrics, and densities on the two surfaces are fixed once and for
all. These auxiliary choices are arbitrary, and formal self-adjointness of $P_b$ need not hold. The Cauchy problem gives bounded isomorphisms
$U_b:\calH^r(\Sigma_-)\to\calH^r(\Sigma_+)$, which are locally uniformly bounded in
$b$, and $b\mapsto U_bf$ is continuous for each $f\in\calH^r(\Sigma_-)$.
These assertions follow from the usual energy estimates (see, for example, \cite{TaylorPsiDO}) and density of
smooth data. 

For fixed $b$, Cauchy evolution admits a microlocal splitting that encodes the propagation of singularities. Namely,
$U_b$ is, modulo a smoothing operator, the sum
$U_b^{\mathrm{pos}} + U_b^{\mathrm{neg}}$ of two Fourier integral
operators (FIOs) whose canonical relations are graphs of canonical maps. We denote these relations by
$\mathcal C_b^\sigma\subset \dotT\Sigma_+\times\dotT\Sigma_-$, where we write as usual $\dotT\Sigma_\pm=T^*\Sigma_\pm \setminus0$. They have an explicit geometric description as follows.

\begin{figure}[ht]
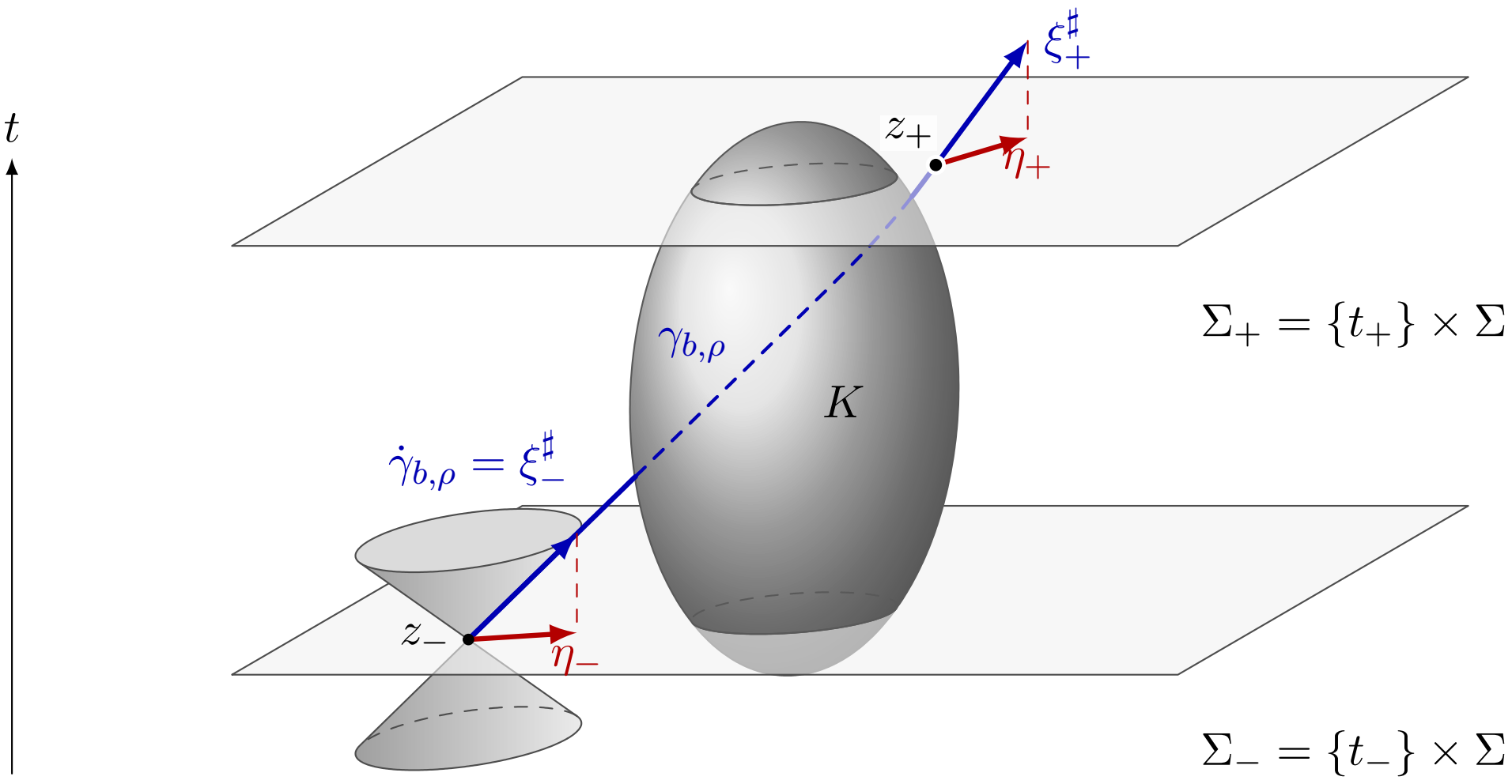

\centering

\caption{The slab between $\Sigma_-$ and $\Sigma_+$ with compact
perturbation support $K$. The lifted null bicharacteristic carries
$\xi_-$ at $z_-$ to $\xi_+$ at $z_+$. The arrows labelled $\eta_\pm$ represent the
restricted covectors $\eta_\pm=\xi_\pm|_{T_{z_\pm}\Sigma_\pm}$ via
$g_b|_{T\Sigma_\pm}$; the dashed segments indicate the tangential
projections of $\xi_\pm^\sharp$ to $T\Sigma_\pm$. Thus
$\rho_\pm=(z_\pm,\eta_\pm)$ give a pair
$(\rho_+,\rho_-)\in\mathcal C^{\mathrm{neg}}$, the Cauchy propagation
relation for the future-directed sheet (see
Section~\ref{sec:reconstruction}). Here $\rho=\rho_-$, and the
projected Hamiltonian orbit $\gamma_{b,\rho}$ defines $A_{b,\rho}$
in \eqref{eq:introA}.}
\label{fig:slab}
\end{figure}

For $\rho_-=(y_-,\eta_-)\in\dotT\Sigma_-$, lift $\eta_-$ to the unique future-directed
null covector $\xi_-$ whose restriction to $T_{y_-}\Sigma_-$ is $\eta_-$.
The Hamiltonian flow of \eqref{eq:H} from this lift projects to a
future-directed, affinely parametrised null geodesic $\gamma_{b,\rho}$ emanating from $\Sigma_-$ and ending on $\Sigma_+$ with a null covector $\xi_+$ (see Fig. \ref{fig:slab}). 
We restrict the transported covector to the final surface to obtain $\rho_+ = (y_+,\eta_+) \in \dotT\Sigma_+$.
The map $\rho_- \mapsto \rho_+$ is canonical, and its graph $\mathcal C_b^{\mathrm{neg}}$ consists of the pairs
$(\rho_+,\rho_-) \in \dotT\Sigma_+\times\dotT\Sigma_-$ obtained in this process.
Similarly, lifting to a past-directed null covector and following its null bicharacteristic flow to $\Sigma_+$ gives the canonical relation $\mathcal C_b^{\mathrm{pos}}$. We call the union
$\mathcal C_b=\mathcal C_b^{\mathrm{pos}}\cup
\mathcal C_b^{\mathrm{neg}}$ the \emph{Cauchy propagation relation}, and its two components the two sheets.
The labels \(\mathrm{pos}\) and \(\mathrm{neg}\) are chosen to agree with the signs of energy in the ultrastatic case (see Section~\ref{sec:ends}). With our metric signature and Fourier convention, the ultrastatic evolution \(T_\pm(t)=\exp(-\rmi t B_\pm)\) associates positive spectral energy with past-directed characteristic covectors.
This FIO description is standard; see \cites{MR0388464,MR1362544}. In Section~\ref{sec:kernel}, we establish its parameter-dependent extension; see in particular \eqref{eq:Cauchyrelation}--\eqref{eq:UWF}, which extend the construction in \cite{StrohmaierStress}*{Appendix C}.
Note that $U_b f$ does not in general depend smoothly on $b$ if one regards it as a function in a fixed Sobolev space. Differentiating in $b$ typically causes a loss of Sobolev regularity.

We now define a smooth map from  $\dotT\Sigma_-$ to the space of one-forms on $B$ by
\[
 A:  \dotT\Sigma_- \to \Omega^1(B),\quad
       \rho \mapsto A_\rho = (b \mapsto A_{b,\rho}),
\]
where
\begin{equation}\label{eq:introA}
 A_{b,\rho}(v)=\frac12\int_{\gamma_{b,\rho}}
 h_b(v)(\dot\gamma_{b,\rho},\dot\gamma_{b,\rho})\,\der r,
 \qquad v\in T_bB.
\end{equation}
Note that $A_{b,\rho}$ is a covector in $T_b^*B$: it is linear in the
parameter tangent vector $v$. The affine scale is fixed by the initial
null covector. Replacing $\rho$ by $c\rho$, for $c>0$, multiplies
$A_{b,\rho}$ by $c$.

The future-directed lift defining $\gamma_{b,\rho}$ belongs to the negative-energy sheet and has parameter covector $A_{b,\rho}$. Reversing its covector and the Hamiltonian parameter gives the positive-energy sheet, with initial spatial covector $-\rho$ and parameter covector $-A_{b,\rho}$. Define the conic sets
\begin{equation}\label{eq:introGamma}
 \Gamma_{\mathrm{pos}}
 =\{(b,-A_{b,\rho}):b\in B,\ \rho\in\dotT\Sigma_-\},
 \qquad \Gamma_{\mathrm{neg}}=-\Gamma_{\mathrm{pos}},
\end{equation}
where the minus sign reverses covectors in each fibre of $T^*B$. These projected sets need not be disjoint.

The following assumption ensures that the parameter space is sufficiently large for our analysis.
\begin{assumption}\label{ass:A}
For every $\rho\in\dotT\Sigma_-$, the one-form $A_\rho$ is nowhere vanishing:
\begin{equation}\label{eq:introNonzero}
 A_{b,\rho}\ne0\quad\text{for every }b\in B
                  \text{ and }\rho\in\dotT\Sigma_-.
\end{equation}
\end{assumption}
This means that for each ray some metric variation in the family has a
nonzero integral; the variation may depend on the ray. Under this
hypothesis the sets in \eqref{eq:introGamma} are closed in
$\dotT B=T^*B\setminus0$. Assumption A can always be achieved by adding a single parameter, provided variations near an intermediate Cauchy surface are allowed; see Section~\ref{sec:functoriality}.

We denote the
kernel-valued distribution associated with $U$ by $K_U$. It takes values in distributional sections of
\[
 \calF=\Hom(\pi_-^*(E_-\oplus E_-),\pi_+^*(E_+\oplus E_+))
       \longrightarrow\Sigma_+\times\Sigma_-,
\]
where $\pi_\pm$ are the two projections and fixed densities are understood. The theorem below shows that parameter testing gives smooth sections of this kernel bundle.

\begin{theorem}[Exact parameter wavefront sets]\label{thm:main}
Suppose Assumption~\ref{ass:A} holds. Then, for every $r,s\in\R$, we have
\begin{gather*}
 K_U\in\calD'(B;C^\infty(\Sigma_+\times\Sigma_-;\calF)),\\
 U\in\calD'(B;\calL(\calH^r(\Sigma_-),\calH^s(\Sigma_+))),\\
 \WF_{\calL(\calH^r,\calH^s)}(U)
 =\Gamma_{\mathrm{pos}}\cup\Gamma_{\mathrm{neg}}.
\end{gather*}

Let $\Pi_b^{\mathrm{pos}}$ and $\Pi_b^{\mathrm{neg}}$ be smooth families of
microlocal Cauchy projections selecting the respective null sheets to all
orders. These projections are idempotent, their mixed products vanish,
and their sum is the identity, all modulo smoothing operators whose
kernels are jointly smooth in $b$ and the spatial variables, as described
in Section~\ref{sec:kernel}.
Put $U_b^\sigma=U_b\Pi_b^\sigma$. Then
$U_b=U_b^{\mathrm{pos}}+U_b^{\mathrm{neg}}$ modulo a jointly smooth kernel, and
\[
 \WF_{\calL(\calH^r,\calH^s)}(U^\sigma)=\Gamma_\sigma,
 \qquad\sigma\in\{\mathrm{pos},\mathrm{neg}\}.
\]
\end{theorem}

The projections are described in Sections~\ref{sec:ends} and~\ref{sec:kernel}.

The first-order Dirac evolution satisfies the same conclusion. Its Cauchy data consist of a single section rather than a pair. We denote this evolution by $V_b$ and its two components, obtained by selecting the null sheets as in Section~\ref{sec:Dirac}, by $V_b^\sigma$.

\begin{corollary}\label{cor:introDirac}
Let $D_b$ be a smooth family of Dirac-type operators on $E$ such that $P_b=D_b^2$ has principal symbol $-g_b^{-1}\Id_E$. Under Assumption~\ref{ass:A}, for every $r,s\in\R$,
\begin{gather*}
 V\in\calD'(B;\calL(H^r(\Sigma_-;E_-),H^s(\Sigma_+;E_+))),\\
 \WF_{\calL(H^r,H^s)}(V)=\Gamma_{\mathrm{pos}}\cup\Gamma_{\mathrm{neg}},
 \qquad
 \WF_{\calL(H^r,H^s)}(V^\sigma)=\Gamma_\sigma.
\end{gather*}
\end{corollary}

The theorem says, in particular, that
$\int\varphi(b)U_b\,\der b$ has a smooth kernel for every
$\varphi\in C^\infty_0(B)$. This integral initially acts on fixed Sobolev
data by strong integration; the resulting operator extends to a smoothing
operator. It is worth noting that $U_b$ itself is neither smoothing nor smooth in $b$ as an operator between fixed Sobolev spaces.

The ordinary Cauchy propagation relation can also be recovered from the metric dependence. Indeed, for any pair of nonzero endpoint covectors $\rho_\pm^0=(z_\pm,\eta_\pm^0)$ and a fixed characteristic sheet $\sigma$, Theorem~\ref{prop:reconstruction} constructs an explicit three-parameter family of diffeomorphism-induced metric variations, chosen for this pair, such that, with $a_\pm=|\eta_\pm^0|$,
\[
 (\rho_+^0,\rho_-^0)\in\mathcal C^\sigma
 \quad\Longleftrightarrow\quad
 (0;a_+,a_++a_-,a_-)\in\WF_{\calL(\calH^r,\calH^s)}(U^\sigma).
\]
Here $\mathcal C^\sigma$ is the ordinary Cauchy propagation relation on sheet $\sigma$. 

In the stress-energy framework motivating this work~\cite{StrohmaierStress}, the stress-energy tensor generates the infinitesimal response to metric variations. The ray functionals here describe the relation of such variations to null propagation. Positivity of these integrals leads to an averaged null-positive condition on metric variations and to polar-cone bounds on the parameter wavefront sets, analogous in form to spectrum conditions. To express this, put
\begin{align}\label{eq:introCone}
 C_b&=\{v\in T_bB:A_{b,\rho}(v)\ge0
                         \text{ for every }\rho\in\dotT\Sigma_-\},\notag\\
 C_b^\circ&=\{\tau\in T_b^*B:\tau(v)\le0
                                      \text{ for every }v\in C_b\}.
\end{align}
Thus $C_b$ is the cone of averaged null-positive metric variations and
$C_b^\circ$ is its polar with the displayed sign convention. Write
$C^\circ=\bigcup_{b\in B}C_b^\circ$. 

\begin{corollary}[Averaged null-positive polar-cone bound]\label{cor:introPolar}
Suppose Assumption~\ref{ass:A} holds. Then, for any $r,s \in \R$
\begin{align*}
 \WF_{\calL(\calH^r,\calH^s)}(U^{\mathrm{pos}})&\subset C^\circ,\\
 \WF_{\calL(\calH^r,\calH^s)}(U^{\mathrm{neg}})&\subset -C^\circ,\\
 \WF_{\calL(\calH^r,\calH^s)}(U)&\subset C^\circ\cup(-C^\circ).
\end{align*}
\end{corollary}

The same polar-cone inclusions hold for the Dirac evolution and its two components in Corollary~\ref{cor:introDirac}.

\begin{remark}
Corollary~\ref{cor:introPolar} also has an interpretation
in quantum field theory on curved spacetimes.
Since the polar cone is convex, the corresponding
wavefront set inclusion is preserved under finite
tensor products. Under the usual implementability
and smoothness assumptions, this yields the same
polar-cone bound for matrix elements of the
second-quantised time evolution between
finite-particle vectors.
\end{remark}

Appendix~C of the earlier stress-energy work \cite{StrohmaierStress} constructs the scalar parameter-dependent canonical relation, its joint Lagrangian, and the fundamental solutions through a parametrix and smoothing correction. Here we use that machinery to obtain exact operator-valued parameter wavefront sets, reconstruction of the Cauchy propagation relation, and compatibility under noncharacteristic restrictions. We include a self-contained proof of the parameter-dependent construction of fundamental solutions for normally hyperbolic operators on vector bundles in Appendix~\ref{app:kernel}, with joint control of parameter dependence, because the exact operator-valued wavefront equalities rely on the nonvanishing of the corresponding bundle-valued principal symbols.

The paper is organised as follows. Section~\ref{sec:ends} discusses
ultrastatic ends, where metric variations recover the time-translation
generators and the energy splitting; it serves as motivation.
Section~\ref{sec:kernel} describes the parameter-dependent fundamental
solutions and the resulting decomposition of the Cauchy evolution.
In Section~\ref{sec:sobolev} we pass to operator-valued wavefront sets,
prove Theorem~\ref{thm:main} and Corollary~\ref{cor:introPolar}, and show
by an example that Assumption~\ref{ass:A} cannot be dropped.
Section~\ref{sec:functoriality} treats distributions on the space of metrics.
Section~\ref{sec:reconstruction} contains the reconstruction of the Cauchy
propagation relation, and Section~\ref{sec:Dirac} the Dirac evolution.
The construction of the fundamental solutions is given in
Appendix~\ref{app:kernel}.

\subsection{Notations and Conventions}

Throughout the symbol $A \subset B$ means that $A$ is a subset of $B$ including possibly the case $A=B$. Complex Hilbert space inner products are conjugate linear in the first argument and linear in the second.
We use the signature convention $(+1,-1,\ldots,-1)$.  The d'Alembert operator is then given by
$\Box = - \delta \der$, in local coordinates $\frac{1}{\sqrt{|g|}} \partial_j(g^{jk}\sqrt{|g|} \partial_k)$.  Spaces of functions such as $C^\infty(M)$ are complex valued unless otherwise stated, e.g. $C^\infty(M)=C^\infty(M,\C)$, whereas the space of real valued smooth functions is denoted by $C^\infty(M,\R)$. The space of distributions on $M$ associated with the test function space $C^\infty_0(M)$ will be denoted by $\mathcal{D}'(M)$. The space of distributions associated with the test function space $C^\infty(M)$ is the space of compactly supported distributions $\mathcal{E}'(M) \subset \mathcal{D}'(M)$.
When it comes to wavefront sets it is convenient to introduce the notation $\dot T^*M = T^*M \setminus 0$ for the cotangent bundle with its zero section removed.
We will also use the notation $\mathrm{I}^*(M,\Lambda)$ for Lagrangian distributions and $\mathrm{I}^*(M \times N, \Lambda')$ for Fourier integral operators from $N$ to $M$ associated to the canonical relation $\Lambda$ as in \cite{MR2512677} to which we refer for details. The notation $\Lambda' = \{(x,\xi,y,-\eta) \mid (x,\xi,y,\eta) \in \Lambda\}$ will be used throughout for canonical relations as well as for wavefront sets of distributional kernels. In particular if $K \in \mathrm{I}^*(M \times N, \Lambda')$ then $\mathrm{WF}(K) \subset \Lambda'$ whereas
$\mathrm{WF}'(K) \subset \Lambda$. For a vector $\xi \in \R^d$ the Japanese bracket $\langle \xi \rangle$ denotes $\langle \xi \rangle = (1+ |\xi|^2)^\frac{1}{2}$ and we also write $\langle (\xi_1, \ldots, \xi_n) \rangle = (1+ |\xi_1|^2 + \ldots + |\xi_n|^2)^\frac{1}{2}$ if we have several vectors, $\xi_j \in \R^{d_j}$.

\section{Time translations and polarisation from metric variations}\label{sec:ends}

We first explain how the usual energy generators are included among metric
variations. For the time-translation example write $g=g_0$ for the background metric and take $P=\Box_g+m^2$ with $m>0$.
The operator $P$ is natural: for diffeomorphism-induced metric families it is transported by pullback. Suppose that the background is ultrastatic at both ends:
\begin{equation}\label{eq:static}
 g=\der t^2-h_-\quad(t\ll0),\qquad
 g=\der t^2-h_+\quad(t\gg0).
\end{equation}
In this scalar example we take the auxiliary connection to be the trivial
connection on the ultrastatic ends, so the Cauchy data there are
$(u,\partial_tu)$, and use the Riemannian densities of $h_\pm$ for the
energy inner products. The connection is extended arbitrarily through the
interior. The assumption $m>0$ is made only for presentational simplicity to avoid zero modes.

Choose smooth functions $0\le\chi_\pm(t)\le1$ with disjoint supports,
with $\chi_+=1$ sufficiently far in the future and $\chi_-=1$ sufficiently
far in the past. Put
\[
 X_\pm=\chi_\pm(t)\partial_t
\]
and denote their flows by $\Phi_\pm^s$. The vector fields are complete
since the functions $\chi_\pm$ are bounded. Their disjoint supports imply
$[X_+,X_-]=0$, so the flows commute.
Write $b=(b_1,b_2)$ for the two displacement parameters. The formula
\[
 F_{b_1,b_2}=\Phi_+^{b_1}\circ\Phi_-^{b_2}
\]
defines a smooth action of $\R^2$. Set $g_{b_1,b_2}=F_{b_1,b_2}^*g$. Near the
future end $F_{b_1,b_2}$ is translation by $b_1$, and near the past end it is
translation by $b_2$. Both $\mathcal L_{X_+}g$ and $\mathcal L_{X_-}g$ have
compact support, since time translation is an isometry of each end.
Consequently, for parameters in a sufficiently small neighbourhood of
zero, $g_{b_1,b_2}-g$ is supported in a fixed compact slab. After shrinking
the parameter neighbourhood, these metrics satisfy the assumptions of
Section~\ref{sec:setup}. Moreover,
\[
 \partial_{b_1} g_{b_1,b_2}=\mathcal L_{X_+}g_{b_1,b_2},\qquad
 \partial_{b_2} g_{b_1,b_2}=\mathcal L_{X_-}g_{b_1,b_2}.
\]

Write $T_\pm(c)$ for Cauchy evolution by time $c$ in the two ultrastatic
ends, and let $U=U_g$. Choose the fixed Cauchy surfaces sufficiently far
inside these ends. Pullback of solutions gives
\begin{equation}\label{eq:conjugation}
 U_{g_{b_1,b_2}}=T_+(b_1)\,U\,T_-(-b_2).
\end{equation}
Indeed, the initial data of $F_{b_1,b_2}^*u$ are $T_-(b_2)$ times those of $u$,
whereas its final data are $T_+(b_1)$ times those of $u$. This proves
\eqref{eq:conjugation} with the indicated sign at the initial end.

In particular, both end evolutions can be recovered from metric variations
and the fixed comparison map $U$:
\begin{equation}\label{eq:recover}
 U^{-1}U_{g_{0,b_2}}=T_-(-b_2),\qquad
 U_{g_{b_1,0}}U^{-1}=T_+(b_1).
\end{equation}

To identify their generators, put
\[
 \mathfrak{H}_\pm=(-\Delta_{h_\pm}+m^2)^{1/2},
\]
where $-\Delta_{h_\pm}$ is nonnegative. On the energy space
$H^1(\Sigma_\pm)\oplus L^2(\Sigma_\pm)$, with squared norm
$\|\mathfrak{H}_\pm u\|^2+\|v\|^2$, the end evolution is the unitary group
\[
 T_\pm(c)=\exp(-\rmi c B_\pm),\qquad
 B_\pm=\rmi \begin{pmatrix}0&\Id\\-\mathfrak{H}_\pm^2&0\end{pmatrix}.
\]
Here $B_\pm$ is self-adjoint on its natural domain. The generators are recovered by
differentiation in $b_1$ and $b_2$ in \eqref{eq:recover}, initially on smooth data. In particular,
\[
 B_-=\frac1{\rmi}\left.\frac{\der}{\der b_2}
       (U^{-1}U_{g_{0,b_2}})\right|_{b_2=0},\qquad
 B_+=\rmi\left.\frac{\der}{\der b_1}
       (U_{g_{b_1,0}}U^{-1})\right|_{b_1=0}.
\]

The spectral projections of these generators give the energy splitting.
At either end, suppressing its subscript, they are
\begin{equation}\label{eq:Pi}
 \Pi^{\mathrm{pos}}=\frac12
 \begin{pmatrix}\Id&\rmi\, \mathfrak{H}^{-1}\\-\rmi\, \mathfrak{H}&\Id\end{pmatrix},
 \qquad
 \Pi^{\mathrm{neg}}=\frac12
 \begin{pmatrix}\Id&-\rmi \,\mathfrak{H}^{-1}\\\rmi \,\mathfrak{H}&\Id\end{pmatrix}.
\end{equation}
Thus the splitting can be defined by metric variations alone: first recover
the generators from \eqref{eq:recover}, then take their positive and negative
spectral projections. Formula~\eqref{eq:Pi} identifies these with the usual
microlocal Cauchy projections. They are bounded on every $\calH^r$ and have
classical pseudodifferential entries of orders $0,-1,1,0$. 

In general, if $X$ generates
a family of diffeomorphisms, its metric variation is $\mathcal L_Xg$. Along
an affinely parametrised null geodesic with tangent $k$,
\begin{equation}\label{eq:boundary}
 \frac12(\mathcal L_Xg)(k,k)=\frac{\der}{\der r}g(X,k),\qquad
 \frac12\int_{r_-}^{r_+}(\mathcal L_Xg)(k,k)\,\der r
 =g(X,k)|_{r_-}^{r_+}.
\end{equation}
For the Hamiltonian \eqref{eq:H}, $k=\xi^\sharp$. In the integral above, $r_-$ and $r_+$ are the Hamiltonian parameters at the initial and final surfaces; on the positive-energy sheet $r_+<r_-$. The parameter-covector
formula proved in Section~\ref{sec:kernel} consequently gives, for the
family $g_{b_1,b_2}$,
\begin{equation}\label{eq:endcov}
 \tau_{b_1}=\xi_+(\partial_t),\qquad
 \tau_{b_2}=-\xi_-(\partial_t).
\end{equation}
These are the temporal components of the final and initial kernel
covectors, with the usual twist in the initial variable. 

The two end translations thus recover the temporal components of the
kernel covectors. General endpoint variations can be used to recover the full relation,
as shown in Section~\ref{sec:reconstruction}.

\section{Parameter-dependent fundamental solutions}\label{sec:kernel}

We recall the joint kernel construction and its geometric description. For each
$b$, let $\Phi_b^r$ be the Hamiltonian flow of $H_b$ on $\dotT M$. Write
$(x,\xi)\sim_b(x',\xi')$ when the two nonzero null covectors lie on the same
orbit. The corresponding relation is Lagrangian for the difference of the
two cotangent symplectic forms.

Choose a common Cauchy surface $\Sigma$. For
$\rho=(y,\eta)\in\dotT\Sigma$, let $\eta_b^\sigma$ be its unique null lift
on one of the two sheets, and set
\[
 \rho_b^\sigma(r)=\Phi_b^r(y,\eta_b^\sigma),\qquad
 \gamma_b^\sigma(r)=\pi_M\rho_b^\sigma(r).
\]
Here $\sigma=\mathrm{pos}$ denotes the past-directed covector sheet, and
$\sigma=\mathrm{neg}$ the future-directed sheet. For increasing $r$, the projected flow is past-directed on the positive-energy sheet and future-directed on the negative-energy sheet. The flow parameters are
restricted to their maximal domains; null-geodesic completeness is not
assumed. Define the covector $\tau_b^\sigma\in T_b^*B$ by
\begin{align}\label{eq:tau}
 \tau_b^\sigma(v;\rho,r_1,r_2)
 &=-\int_{r_2}^{r_1}\der_bH_b(v)(\rho_b^\sigma(r))\,\der r\notag\\
 &=\frac12\int_{r_2}^{r_1}
 h_b(v)(\dot\gamma_b^\sigma(r),\dot\gamma_b^\sigma(r))\,\der r.
\end{align}
Differentiating the inverse metric at fixed $(x,\xi)$ gives
$\der_bH_b(v)=-\tfrac12h_b(v)(\xi^\sharp,\xi^\sharp)$, and
$\dot\gamma=\xi^\sharp$.

Let $\Lambda$ consist of the points
\[
 (b,\tau_b^\sigma;\rho_b^\sigma(r_1);\rho_b^\sigma(r_2)),
\]
for both sheets. It is understood as a canonical relation, with a minus
sign on the symplectic form in the final displayed factor. Its twist
\[
 \Lambda'=\{(b,\tau;x,\xi;x',-\xi'):
                   (b,\tau;x,\xi;x',\xi')\in\Lambda\}
\]
is a homogeneous Lagrangian in $\dotT(B\times M\times M)$.

To see this, regard $H(b,\tau,x,\xi)=H_b(x,\xi)$ as a Hamiltonian on
$T^*(B\times M)$. Its Hamiltonian vector field is
\[
 X_H=X_{H_b}-\sum_j(\partial_{b^j}H_b)\partial_{\tau_j}.
\]
Starting from the null part of the identity relation with $\tau=0$, the
flow-out therefore has parameter component \eqref{eq:tau}. The initial
submanifold is isotropic and the flow-out has the required Lagrangian
dimension. Homogeneity follows by scaling the spacetime covectors by
$c>0$ and the affine flow intervals by $c^{-1}$; the parameter covector then
scales by $c$. Smooth dependence follows from the smooth dependence of the
Hamiltonian flow on $b$.

Set
\[
 \Delta^*=\{(b,0;x,\xi;x,-\xi):(x,\xi)\in\dotT M\}.
\]
The kernels are distributional sections of $E\boxtimes E^*$, whose fibre at $(x,x')$ is $\Hom(E_{x'},E_x)$; pullback to the product with $B$ is understood. We identify kernels with distributions using the fixed density of $g_0$.
Changing that density multiplies kernels by smooth nonvanishing factors.

\begin{theorem}\label{thm:kernel}
The retarded and advanced fundamental solutions have joint kernels in the space of distributions
$\calD'(B\times M\times M;E\boxtimes E^*)$ and satisfy
\[
 \WF(G_{b,\mathrm{ret/adv}})\subset\Lambda'\cup\Delta^*.
\]
The Pauli--Jordan propagator defined by $G_b=G_{b,\mathrm{ret}}-G_{b,\mathrm{adv}}$ belongs to the space of Lagrangian distributions
\[
 I^{-\frac32-\frac{\dim B}{4}}(B\times M\times M,\Lambda';E\boxtimes E^*)
\]
with a classical polyhomogeneous symbol.
\end{theorem}

The proof of Theorem~\ref{thm:kernel} is given in Appendix~\ref{app:kernel}.

Fix a relatively compact parameter neighbourhood $B_0\Subset B$ and a
finite time interval $I=(a,c)$, with all estimates taken on compact
subintervals and compact parameter subsets. Parallel
transport for the fixed connection along the curves $t\mapsto(t,x)$
identifies $E|_{I\times\Sigma}$ with the pullback of a bundle on $\Sigma$.
This identification is independent of $b$. Choose a smooth positive
density $\der\mu$ on $\Sigma$ and use the product density
$\der t\,\der\mu$ in this proof. Changing back to the fixed spacetime
density only multiplies the kernels by a fixed smooth nonvanishing factor.

In these identifications we may write
\begin{equation}\label{eq:timeform}
 P_b=c_b(t)\partial_t^2+A_{b,1}(t)\partial_t+A_{b,2}(t),
\end{equation}
where $c_b$ is a nowhere vanishing scalar multiplication operator and
$A_{b,j}$ is a tangential differential operator of order at most $j$.
All coefficients are smooth in $(b,t,x)$. The fact that $c_b$ is invertible
follows from the assumption that $\der t$ is timelike. Let $\delta_\Sigma$
denote the kernel of the identity on sections over $\Sigma$, relative to
$\der\mu$. Equal-time identities below are identities of distribution
kernels on $\Sigma\times\Sigma$.

\subsection{Cauchy evolution}

In the product-density
convention above, let $G_b(t,s)$ denote the Pauli--Jordan kernel $G_b$ restricted to
the indicated times and viewed as an operator on $\Sigma$. The solution
with $u|_{t=s}=f$ and $\partial_tu|_{t=s}=g$ is
\begin{equation}\label{eq:CauchyfromGreen}
 u(t)=-\partial_sG_b(t,s)c_b(s)f
       +G_b(t,s)\big(c_b(s)g+
                    (A_{b,1}(s)-\partial_sc_b(s))f\big).
\end{equation}
For verification, extend a homogeneous solution by zero to the past of
$s$. Applying~\eqref{eq:timeform} gives the distributional source
\[
 \delta'(t-s)c_b(s)f+
 \delta(t-s)\big(c_b(s)g+
                    (A_{b,1}(s)-\partial_sc_b(s))f\big).
\]
The retarded solution of this source is $\theta(t-s)u$, the advanced
solution is $-\theta(s-t)u$, and their difference is $u$.
Applying $G_b$ to the source gives~\eqref{eq:CauchyfromGreen}.
The conversion from $\partial_tu$ to $\nabla_{n_b}u$ uses only tangential
operators and multiplications with smoothly parameter-dependent
coefficients. The source restrictions and derivatives in
\eqref{eq:CauchyfromGreen} have no zero input spatial covectors in their
joint wavefront relations. Consequently smooth Cauchy data depending
smoothly on $b$, and on any additional variables, give jointly smooth
solutions. This conclusion also follows from the Green kernel construction with parameters.

Let $\mathcal U_b=R_{b,\Sigma_-}^{-1}$ be the solution operator for initial
Cauchy data. Appendix~\ref{app:kernel} gives a decomposition
$G_b\equiv G_b^{\mathrm{pos}}+G_b^{\mathrm{neg}}$ on the two null sheets,
where $\equiv$ denotes equality modulo jointly smooth kernels. Apply
\eqref{eq:CauchyfromGreen} to each $G_b^\sigma$ with $s=t_-$, using the
conversion of initial data above, and call the resulting operators
$\mathcal U_b^\sigma$. They sum to $\mathcal U_b$ modulo a jointly smooth
kernel. The Cauchy restrictions and normal derivatives are well-defined
since a nonzero null covector cannot be conormal to a spacelike hypersurface
(\cite{Ho1}*{Theorem 8.2.4}). Set
\[
 \Pi_b^\sigma=R_{b,\Sigma_-}\mathcal U_b^\sigma.
\]
At equal times both sheets restrict to the spatial diagonal with zero
parameter covector. Thus the $\Pi_b^\sigma$ are smooth families of
pseudodifferential matrices, and $R_{b,\Sigma_-}\mathcal U_b=\Id$ gives
$\Pi_b^{\mathrm{pos}}+\Pi_b^{\mathrm{neg}}\equiv\Id$.

Since $P_b\mathcal U_b=0$ and the two spacetime sheets are disjoint,
$P_b\mathcal U_b^\sigma$ has a jointly smooth kernel. The difference
$\mathcal U_b^\sigma-\mathcal U_b\Pi_b^\sigma$ has zero initial data and
jointly smooth forcing. Its kernel is jointly smooth by the inhomogeneous
Cauchy formula, which integrates $G_b(t,r)$ against the forcing from $s$
to $t$. The integrand is jointly smooth by the input wavefront condition
above. Writing $\bar\sigma$ for the other sheet, we obtain
\[
 \mathcal U_b^\sigma-\mathcal U_b^\sigma\Pi_b^\sigma
 \equiv \mathcal U_b^{\bar\sigma}\Pi_b^\sigma.
\]
Composition on the right with $\Pi_b^\sigma$ preserves the output sheet.
The two sides thus have disjoint wavefront sets and are jointly smooth.
Taking their initial Cauchy traces gives
$(\Pi_b^\sigma)^2\equiv\Pi_b^\sigma$ and
$\Pi_b^{\bar\sigma}\Pi_b^\sigma\equiv0$.
On each sheet the normal trace multiplies the principal symbol of the
value by $\rmi\xi(n_b)$. The symbols of the two initial traces therefore
take values in the respective Cauchy polarisations and, since they sum
to the identity, are the projections onto them. In the ultrastatic
setting we may take the exact projections~\eqref{eq:Pi}, which select
the same sheets.

We now take the Cauchy data at the final surface. Write
$U_b^\sigma=U_b\Pi_b^\sigma\equiv R_{b,\Sigma_+}\mathcal U_b^\sigma$.
The canonical relation of this restricted branch is
\begin{equation}\label{eq:Cauchyrelation}
 \Lambda_U^\sigma
 =\{(b,\tau_b^\sigma(\rho);\kappa_b^\sigma(\rho);\rho):
                              \rho\in\dotT\Sigma_-\},
\end{equation}
where $\kappa_b^\sigma$ is null propagation to $\Sigma_+$ followed by
restriction of covectors. In \eqref{eq:tau} we now take the initial and
final intersections with these fixed surfaces. Thus $r_1<r_2$ on the positive-energy sheet and $r_1>r_2$ on the negative-energy sheet. Derivatives of their flow
hitting times introduce no additional action terms, since $H_b=0$ on the
null relation.

For Cauchy data $(u,\partial_tu)$, the upper right entry of $U_b^\sigma$
agrees with $G_b^\sigma(t_+,t_-)c_b(t_-)$ modulo a jointly smooth kernel.
The initial normalisation and transport in
Appendix~\ref{app:kernel} give a nonzero symbol for this entry on the
whole relation. The changes of Cauchy data are elliptic isomorphisms,
so
\begin{equation}\label{eq:UWF}
 \WF(K_{U^\sigma})=(\Lambda_U^\sigma)'.
\end{equation}
The two restricted relations may intersect, but their contributions
cannot cancel in the full Cauchy matrix. Indeed, composition on the
right with $\Pi_b^\sigma$ recovers $U_b^\sigma$ modulo a jointly smooth
kernel and does not enlarge the joint wavefront set. Consequently
\[
 \WF(K_U)=(\Lambda_U^{\mathrm{pos}})'\cup(\Lambda_U^{\mathrm{neg}})'.
\]

The Cauchy matrix entries have the usual orders $0,-1,1,0$ at fixed $b$.
Conjugating at either end by
$\operatorname{diag}(Q^{1/2},Q^{-1/2})$, where $Q$ is a fixed positive
elliptic operator of order one on the corresponding surface, gives an
order-zero matrix at fixed $b$ and joint FIO order $-\dim B/4$.
Differentiation in the parameter increases the spatial order by at most
one. Consequently, for fixed data,
\begin{equation}\label{eq:loss}
 f\in\calH^{r+k}(\Sigma_-)
 \quad\Longrightarrow\quad
 b\longmapsto U_bf\in C^k(B;\calH^r(\Sigma_+)).
\end{equation}

\section{Sobolev-valued wavefront sets and averaged null variations}\label{sec:sobolev}

For a Banach space $F$ and $u\in\calD'(B;F)$, we use the norm
definition of the wavefront set. Hence,
$(b_0,\zeta_0)\notin\WF_F(u)$ if and only if there are a cutoff function $\chi \in C^\infty_0(B)$ with $\chi=1$ near $b_0$
and a conic neighbourhood $V$ of $\zeta_0$ such that
\[
 \|\widehat{\chi u}(\zeta)\|_F
 \le C_N\langle\zeta\rangle^{-N},\qquad \zeta\in V,
\]
for every $N$, where
$\widehat{\chi u}(\zeta)=u(\chi(b)e^{-\rmi b\cdot\zeta})$.
Equivalently, one may test against continuous linear functionals $\ell\in F^*$, provided the cutoff $\chi$ and the cone $V$ are the same for all $\ell$. More precisely, for these fixed $\chi,V$ it is enough that, for every $N$ and $\ell$,
\[
 \sup_{\zeta\in V}\langle\zeta\rangle^N
       |\ell(\widehat{\chi u}(\zeta))|<\infty.
\]
Equivalence to the norm-definition follows from the Banach-Steinhaus theorem.
This extends the Hilbert-valued definition; see
\cite{MR1936535}*{Section 2}. For a Fr\'echet space $F$ we require the same
estimate for every continuous seminorm, with a common $\chi$ and $V$.
Every continuous linear map $T:F\to G$ satisfies
$\WF_G(Tu)\subset\WF_F(u)$.

In the following $X,Y$ are compact manifolds with fixed smooth positive densities and finite rank vector bundles $E_X\to X$ and $E_Y\to Y$. Let as usual $E_Y\boxtimes {E_X}^*$ be a vector bundle on $Y \times X$ with fibre $\Hom((E_X)_x,(E_Y)_y)$ at $(y,x)$, so Schwartz kernels of operators from smooth sections of $E_X$ to distributional sections of $E_Y$ are distributions in
$\calD'(Y \times X; E_Y \boxtimes {E_X}^*)$. To avoid heavy notation we will simply write 
$\calD'(B\times Y\times X; E_Y\boxtimes {E_X}^*)$ for the distributional sections of the bundle pulled back to $B\times Y\times X$.

\begin{lemma}\label{lem:kernel}
Let $K\in\calD'(B\times Y\times X; E_Y\boxtimes {E_X}^*)$ be a distributional kernel. Suppose that
\begin{equation}\label{eq:transverse}
 (b,\tau;y,\eta;x,\xi)\in\WF(K)\quad\Longrightarrow\quad\tau\ne0.
\end{equation}
Then integration with a test function defines a distribution with values in smooth kernels,
\[
 \mathcal K:C^\infty_0(B)\longrightarrow C^\infty(Y\times X;E_Y\boxtimes {E_X}^*),\qquad
 \mathcal K(\varphi)=\langle K,\varphi\rangle_B.
\]
The projection
\[
 \Gamma=\{(b,\tau):\exists(y,\eta;x,\xi),\
                         (b,\tau;y,\eta;x,\xi)\in\WF(K)\}
\]
is closed in $\dotT B$.
In the remaining formulas, $H^r(X)$ and $H^s(Y)$ denote Sobolev spaces of sections of $E_X$ and $E_Y$, respectively. In particular, for every $r,s\in\R$, the operators with kernels
$\mathcal K(\varphi)$ define
\[
 \mathcal V\in\calD'(B;\calL(H^r(X),H^s(Y))).
\]
For every $r,s\in\R$, we then have
\begin{equation}\label{eq:kernelWFexact}
 \WF_{C^\infty(Y\times X;E_Y\boxtimes {E_X}^*)}(\mathcal K)
 =\WF_{\calL(H^r,H^s)}(\mathcal V)=\Gamma.
\end{equation}
More precisely, away from $\Gamma$ there are a cutoff $\chi$ and a conic
neighbourhood $V$ for which
\begin{equation}\label{eq:operatorbound}
 \|\mathcal V(\chi(b)e^{-\rmi b\cdot\zeta})\|_{H^r\to H^s}
 \le C_{r,s,N}\langle\zeta\rangle^{-N},\qquad\zeta\in V,
\end{equation}
for every $r,s,N$. If $K$ is
also the joint kernel of a strongly continuous, locally uniformly bounded
family $V_b:H^{r_0}(X)\to H^{s_0}(Y)$, then
\[
 \mathcal V(\varphi)f=\int\varphi(b)V_bf\,\der b,
 \qquad f\in H^{r_0}(X),
\]
where the integral is a strong integral in $H^{s_0}(Y)$. Consequently
$\WF_{H^{s_0}(Y)}(V_\bullet f)\subset\Gamma$, uniformly for $f$ in the
unit ball of $H^{r_0}(X)$.
\end{lemma}

\begin{proof}
We work in local bundle frames. On compact subsets of $B$,
condition~\eqref{eq:transverse} implies
\begin{equation}\label{eq:properprojection}
 |\eta|+|\xi|\le C|\tau|\quad\hbox{on }\WF(K).
\end{equation}
Otherwise, after normalisation and passage to a subsequence, compactness
would give a nonzero element of $\WF(K)$ with $\tau=0$. The same estimate shows
that the projection is closed in $\dotT B$.

For $\varphi\in C^\infty_0(B)$, the pushforward of $\varphi K$ to
$Y\times X$ is well-defined because its support is compact in the
parameter variable. The pushforward wavefront estimate, in the kernel form of \cite{Ho1}*{Theorem 8.2.12}, gives
\[
 \WF(\mathcal K(\varphi))\subset
 \{(y,\eta;x,\xi):\exists b,\
             (b,0;y,\eta;x,\xi)\in\WF(K)\}=\varnothing.
\]
Thus $\mathcal K(\varphi)$ is smooth, and the map
$\varphi\mapsto\mathcal K(\varphi)$ is continuous from $C^\infty_0(B)$
to $C^\infty$ (see \cite{Ho1}*{remark after Theorem 8.2.12}).
In particular, for each $s_1$ and each compact $Q\Subset B$ there are
$k$ and $C$ such that
\[
 \|\mathcal K(\varphi)\|_{H^{s_1}(Y\times X;E_Y\boxtimes {E_X}^*)}
 \le C\max_{|\alpha|\le k}\sup|\partial^\alpha\varphi|,
 \qquad\supp\varphi\subset Q,
\]
in local parameter coordinates. The order $k$ may depend on $s_1$.

Fix $(b_0,\zeta_0)\notin\Gamma$. Choose $\chi=1$ near $b_0$ and a conic
neighbourhood $V$ of $\zeta_0$ so that no element of $\WF(K)$ over
$\supp\chi$ has parameter direction in $\overline V$. Localise the kernel
in finitely many spatial coordinate charts. For each localised kernel
$K_j$, the closed cone of total covectors with
$\zeta\in\overline V\cup\{0\}$ is disjoint from its wavefront directions:
the nonzero $\zeta$ case follows from the choice of $V$, and the $\zeta=0$
case from \eqref{eq:transverse}. The Fourier characterisation of the
ordinary wavefront set, \cite{Ho1}*{Definition 8.1.2 and Lemma 8.1.1}, therefore gives, for every $M$,
\begin{equation}\label{eq:totalFourier}
 |\widehat{\chi K_j}(\zeta,\eta,\xi)|
 \le C_M\langle(\zeta,\eta,\xi)\rangle^{-M},\qquad\zeta\in V,
\end{equation}
uniformly in $\eta,\xi$. Integrating with arbitrary polynomial weights in
the spatial frequencies yields
\[
 \|\mathcal K(\chi(b)e^{-\rmi b\cdot\zeta})\|_{H^{s_1}(Y\times X;E_Y\boxtimes {E_X}^*)}
 \le C_{s_1,N}\langle\zeta\rangle^{-N}
\]
for every $s_1,N$. These estimates prove the smooth-kernel-valued wavefront
inclusion. For sufficiently large $s_1$ the kernel Sobolev norm
controls its operator norm from $H^r(X)$ to $H^s(Y)$. Summing the finite
localisation proves \eqref{eq:operatorbound}.

For the converse, suppose $(b_0,\zeta_0)$ is absent from the operator-valued wavefront set for fixed $r,s$. Choose its cutoff $\chi$ and cone $V$. Each local component of $\widehat{\chi K_j}(\zeta,\eta,\xi)$ is obtained by testing the operator $\mathcal V(\chi e^{-\rmi b\cdot\zeta})$ against spatially cut-off plane waves. Their $H^r$ and $H^{-s}$ norms grow at most polynomially in the spatial frequencies. Thus, for some fixed $k$ and every $N$,
\[
 |\widehat{\chi K_j}(\zeta,\eta,\xi)|
 \le C_N\langle\zeta\rangle^{-N}\langle(\eta,\xi)\rangle^k,
 \qquad\zeta\in V.
\]
Near any joint covector over $b_0$ with parameter component $\zeta_0\ne0$ we may restrict to a cone where $|\eta|+|\xi|\le C|\zeta|$. This estimate then gives rapid decay in the total frequency and excludes that joint covector from $\WF(K)$. Hence $\Gamma\subset\WF_{\calL(H^r,H^s)}(\mathcal V)$. Together with the already proved inclusions and the continuity of the map from smooth kernels to Sobolev operators, this proves~\eqref{eq:kernelWFexact}.

If the strongly continuous
family exists, its strong integral has kernel $\mathcal K(\varphi)$ by
testing against smooth input and output functions; density gives the
identity for all $f\in H^{r_0}(X)$. Applying \eqref{eq:operatorbound} to
$f$ gives the last assertion. 
\end{proof}

\begin{remark}
The operator target in Lemma~\ref{lem:kernel} may also be taken to be the Hilbert space
$\mathcal S_2(H^r(X),H^s(Y))$ of Hilbert--Schmidt operators. Indeed, choose fixed positive elliptic operators $Q_X,Q_Y$ of order one defining the Sobolev norms. For a smooth kernel $k$, the Hilbert--Schmidt norm is
\[
 \|A_k\|_{\mathcal S_2(H^r,H^s)}
   =\|\text{kernel of }Q_Y^sA_kQ_X^{-r}\|_{L^2(Y\times X)}.
\]
For $s_1$ sufficiently large, this is bounded by a constant times
$\|k\|_{H^{s_1}(Y\times X;E_Y\boxtimes {E_X}^*)}$. The lemma and the estimates in its proof therefore give
\[
 \mathcal V\in\calD'(B;\mathcal S_2(H^r(X),H^s(Y))),\qquad
 \WF_{\mathcal S_2(H^r,H^s)}(\mathcal V)=\Gamma.
\]
Since this target is a Hilbert space, its wavefront set is precisely the usual Hilbert-space-valued wavefront set. In the setting of the lemma it agrees with the operator-norm wavefront set. One inclusion follows from the continuous inclusion into $\calL(H^r,H^s)$. For the converse, the operator norm controls a sufficiently weak kernel Sobolev norm: after the fixed elliptic conjugations, compose on the left with $Q_Y^{-s_2}$, where $s_2>\dim Y/2$, and use
\[
 \|Q_Y^{-s_2}A\|_{\mathcal S_2(L^2,L^2)}
       \le\|Q_Y^{-s_2}\|_{\mathcal S_2(L^2,L^2)}\|A\|_{\calL(L^2,L^2)}.
\]
The smooth-kernel-valued distribution property gives polynomial Fourier bounds in every stronger kernel Sobolev norm. Interpolation therefore upgrades rapid decay in the operator norm to rapid decay in $H^{s_1}$, and hence in the Hilbert--Schmidt norm, with the same cutoff and cone.
\end{remark}

The condition in the lemma concerns the projection of the full kernel
relation, not the nonvanishing of its input or output covectors. For the
relation \eqref{eq:Cauchyrelation} it is Assumption~\ref{ass:A}, equivalently
\begin{equation}\label{eq:nonzero}
 \tau_b^\sigma(\rho)\ne0
 \quad\hbox{for all }\rho\in\dotT\Sigma_-.
\end{equation}
Nonvanishing at a fixed parameter point is uniform on a cosphere bundle
and persists on a parameter neighbourhood. The canonical graphs are
homogeneous and have compact base; thus they have the closedness and
compactness properties required by Lemma~\ref{lem:kernel}.

\begin{theorem}\label{thm:Sobolev}
Suppose Assumption~\ref{ass:A} holds on a parameter neighbourhood. Put
\[
 \Gamma_\sigma
 =\{(b,\tau_b^\sigma(\rho)):\rho\in\dotT\Sigma_-\}.
\]
With $F=C^\infty(\Sigma_+\times\Sigma_-;\calF)$ as in
Section~\ref{sec:setup}, the kernel of $U^\sigma$ satisfies
\[
 K_{U^\sigma}\in\calD'(B;F),\qquad
 \WF_F(K_{U^\sigma})=\Gamma_\sigma.
\]
For every $r,s\in\R$, it therefore defines an operator-valued distribution
satisfying
\begin{gather*}
 U^\sigma\in\calD'(B;\calL(\calH^r(\Sigma_-),\calH^s(\Sigma_+))),\\
 \WF_{\calL(\calH^r,\calH^s)}(U^\sigma)=\Gamma_\sigma.
\end{gather*}
For the evolution,
\[
 \WF_F(K_U)=\WF_{\calL(\calH^r,\calH^s)}(U)
 =\Gamma_{\mathrm{pos}}\cup\Gamma_{\mathrm{neg}}.
\]
For every $r\in\R$ and $f\in\calH^r(\Sigma_-)$,
\begin{equation}\label{eq:mainWF}
 \WF_{\calH^r(\Sigma_+)}(U_\bullet^\sigma f)
 \subset\Gamma_\sigma.
\end{equation}
The rapid-decay estimates are uniform over the unit ball of
$\calH^r(\Sigma_-)$. In particular, the same inclusion holds for every
continuous scalar pairing of the evolved data. If both sheets are
included, then
\[
 \WF_{\calH^r(\Sigma_+)}(U_\bullet f)
 \subset\Gamma_{\mathrm{pos}}\cup\Gamma_{\mathrm{neg}},
 \qquad \Gamma_{\mathrm{neg}}=-\Gamma_{\mathrm{pos}}.
\]
\end{theorem}

\begin{proof}
Apply Lemma~\ref{lem:kernel} to \eqref{eq:UWF}, using the fixed elliptic
weights described above. The smoothing part of the initial decomposition
gives a smooth Sobolev-valued function. Reversing the spacetime covectors
and the orientation of the Hamiltonian parameter in \eqref{eq:tau}
gives
\[
 \tau_b^{\mathrm{pos}}(-\rho)=-A_{b,\rho}.
\]
Thus the two projected relations are negatives of one another. The equality for the evolution follows from the joint kernel equality established after~\eqref{eq:UWF} and~\eqref{eq:kernelWFexact}. The assertions for individual initial data follow by evaluation and remain inclusions. The smooth-kernel conclusion and the operator-valued equalities prove Theorem~\ref{thm:main}. Corollary~\ref{cor:introPolar} follows from $-A_{b,\rho}\in C_b^\circ$, which gives $\Gamma_{\mathrm{pos}}\subset C^\circ$.
\end{proof}

One can refine \eqref{eq:mainWF} by retaining only rays whose initial covectors occur in the joint wavefront set of $(b,y)\mapsto\Pi_b^\sigma f(y)$, with the wavefront set of a Cauchy pair defined as the union of those of its components. The parameter covector of this joint distribution is zero, since the projections are smooth families of pseudodifferential operators. When the projections are independent of $b$, this reduces to retaining only rays issuing from $\WF(\Pi^\sigma f)$. To see this, localise the
initial covectors by pseudodifferential cutoffs. The part of the datum
which is microlocally smooth contributes a smooth parameter-dependent
solution, whereas on the remaining compact set of cosphere directions the
proof of Lemma~\ref{lem:kernel} applies to the cut-off evolution. In this
argument it suffices that the nonvanishing condition in Assumption~\ref{ass:A} hold on the rays carrying the
initial singularities.

A convenient sufficient condition for Assumption~\ref{ass:A} near
$b_0$ is the existence of a direction $v_0\in T_{b_0}B$ such that
\begin{equation}\label{eq:strict}
 A_{b_0,\rho}(v_0)>0\quad
                       \text{for every }\rho\in\dotT\Sigma_-.
\end{equation}
By compactness of the unit cosphere bundle and smooth dependence, this
strict inequality persists locally after extending $v_0$ to a smooth
parameter direction. Thus Assumption~\ref{ass:A} holds locally, and
Corollary~\ref{cor:introPolar} applies.

Finite-dimensional cone duality gives
\[
 C_b^\circ=\overline{\operatorname{cone}
                  \{-A_{b,\rho}:\rho\in\dotT\Sigma_-\}},
\]
where $\operatorname{cone}$ denotes finite nonnegative linear combinations.

The ultrastatic construction gives a natural instance of
\eqref{eq:strict}. The final-translation direction $\partial_{b_1}$ in
Section~\ref{sec:ends} satisfies, by \eqref{eq:endcov},
\[
 A_{b,\rho}(\partial_{b_1})=\xi_+(\partial_t)>0
\]
where $\xi_+$ is the future-directed final covector of the ray defining $A_\rho$. Thus adjoining this translation parameter
to any admissible family fixed on the ultrastatic ends supplies the strict
direction.
The conclusion is a wavefront estimate on the enlarged parameter manifold.
Restriction to a smaller parameter manifold obeys the usual
noncharacteristic pullback condition.

\subsection{Failure without Assumption A}\label{sec:counterexample}

The following example on a flat cylinder shows that Assumption \ref{ass:A} is indeed necessary.
On $\R\times S^1$ let
\[
 g_s=(1+s^2\chi(t))^2\der t^2-\der x^2,
 \qquad \chi\in C^\infty_0((0,1)),\quad\chi\ge0,
 \quad\int_0^1\chi(t)\,\der t=1.
\]
This is the slice $a=0$ of the family
\[
 g_{s,a}=(1+(s^2+a)\chi(t))^2\der t^2-\der x^2.
\]
The $a$ direction makes the enlarged family satisfy
Assumption~\ref{ass:A} near $(0,0)$. At $(0,0)$ the $s$ derivative
vanishes, so the ray covectors are multiples of $\der a$, conormal
to the slice.

For the massless wave equation, one polarised evolution from $t=0$ to
$t=1$, after removal of a fixed evolution factor, is $\exp(\rmi s^2D)$
on the Fourier modes $De_n=ne_n$, $n\ge1$. The displacement
$f=\sum_{n\ge1}n^{-1}e_n$ and its polarised normal derivative give data in
$\calH^{-1/2}$. A continuous scalar pairing of the evolved data is
\[
 F(s)=\sum_{n\ge1}\frac{e^{\rmi ns^2}}{n^2}.
\]
For $0<q<2\pi$, the classical Fourier-series identity gives
$\sum_{n\ge1}\cos(nq)/n=-\log(2\sin(q/2))$. Integration from zero yields
\[
 \operatorname{Im}F(s)=-2s^2\log|s|+s^2+O(s^6).
\]
Thus $F$ is not smooth at zero. It is even, so both wavefront directions
occur there. On the other hand $\partial_sg_s|_{s=0}=0$, and every ray
parameter covector vanishes at that point. In particular $C_0=T_0B$ and
$C_0^\circ=\{0\}$.

Parameter averaging need not produce a smoothing operator either. Indeed,
for $\varphi\in C^\infty_0(B)$ with $\varphi(0)\ne0$, the averaged
polarised evolution, after removal of the same fixed evolution factor,
acts on $e_n$ by multiplication with
\[
 a_n=\int_{\R}\varphi(s)e^{\rmi n s^2}\,\der s
     =e^{\rmi\pi/4}\sqrt{\frac{\pi}{n}}\,\varphi(0)
       +O(n^{-3/2}),\qquad n\to\infty,
\]
by one-dimensional stationary phase at $s=0$. Thus the multipliers do not
decay rapidly: averaging gains only half a derivative on these modes and
does not yield a smoothing operator.

\section{Distributions on Fr\'echet manifolds and their wavefront set}
\label{sec:functoriality}

In this section we briefly describe a notion of generalised function on Fr\'echet manifolds that allows to directly treat the Cauchy evolution operator as a distribution on the space of metrics.

Let $Q$ be an open neighbourhood of zero in a Fr\'echet space and let
$V$ be a Banach space. Let $\mathcal M$ be a collection of smooth maps
$f:B_f\to Q$ from finite-dimensional open balls. We will require that this family is sufficiently rich so that it eventually includes all finite dimensional smooth maps as a restriction. To be precise, we require that, given
any smooth map $\tilde f:B\to Q$ from a finite-dimensional manifold
and any $b_0\in B$, there be a neighbourhood $U$ of $b_0$, a map
$f\in\mathcal M$, and a smooth inclusion $i:U\hookrightarrow B_f$ such
that $\tilde f|_U=f\circ i$.

We call $F$ a \emph{tame distribution} on $Q$ with values in $V$ if it
is given by a family $F_f\in\mathcal D'(B_f;V)$, $f\in\mathcal M$,
satisfying
\begin{equation}\label{compatcond}
 \WF_V(F_{f'})\cap N_i=\varnothing,\qquad
 F_f|_U=i^*F_{f'}
\end{equation}
whenever $f,f'\in\mathcal M$, $U\subset B_f$ is open, and
$f|_U=f'\circ i$ for a smooth embedding $i:U\hookrightarrow B_{f'}$.
Here
\[
 N_i=\{(i(b),\eta):b\in U,\ \eta\ne0,\ (di_b)^*\eta=0\}
\]
is the conormal of $i$. The first condition ensures that the pullback
exists.

We write $T_q^*Q$ for the continuous dual of the model space, equipped
with the weak$^*$ topology. For a cone $\mathcal W\subset T^*Q\setminus0$ and
$f\in\mathcal M$ put
\[
 f^*\mathcal W=\{(b,(df_b)^*\xi):b\in B_f,\ \xi\in \mathcal W_{f(b)},\
        (df_b)^*\xi\ne0\}.
\]
We call $\mathcal W$ \emph{admissible} for $F$, and write $\WF_V(F)\subset\mathcal W$, if the fibres $\mathcal W_q$ are weak$^*$-closed and
\[
 \WF_V(F_f)\subset f^*\mathcal W\qquad\text{for every }f\in\mathcal M.
\]
If there is a smallest admissible cone, we call it the wavefront set
$\WF_V(F)$. If $Q$ is finite-dimensional and $F_f=f^*F$, this is the
usual wavefront set. Indeed, the usual wavefront set satisfies the
condition by the pullback theorem \cite{Ho1}*{Theorem 8.2.4}, whose
proof carries over to Banach-space-valued distributions. Conversely, the
richness of $\mathcal M$ provides maps $f\in\mathcal M$ that are
submersions near any given point. For these $(df_b)^*$ is injective and
$\WF_V(F_f)=f^*\WF_V(F)$, so every admissible cone contains the usual
wavefront set. 

For the Cauchy evolution, take $Q$ to be a local chart of metrics equal to
$g_0$ off $K$, with $K$ enlarged as described below. We use a common
Cauchy temporal function $t$ and fix a smooth assignment $g\mapsto P_g$
with $P_g=P_0$ off $K$. For fixed $r,s$ take
$V=\calL(\calH^r(\Sigma_-),\calH^s(\Sigma_+))$, and let $\mathcal M$ consist
of the maps $f:B_f\to Q$ satisfying Assumption~\ref{ass:A}. Write $U_f$
for the Cauchy evolution on such a family, regarded as a distribution
with values in $V$ by Theorem~\ref{thm:Sobolev}.

If $f=f'\circ i$ for two such maps, the chain rule in~\eqref{eq:introA}
gives
\[
 A^f_{b,\rho}=(di_b)^*A^{f'}_{i(b),\rho}.
\]
The left-hand side is nonzero by Assumption~\ref{ass:A} for $f$.
The exact wavefront formula for $U_{f'}$ therefore gives
$\WF_V(U_{f'})\cap N_i=\varnothing$. The continuous injection of $V$
into $\calD'(\Sigma_+\times\Sigma_-;\calF)$ commutes with pullback.
The joint kernels have no purely parameter covectors
by~\eqref{eq:Cauchyrelation} and~\eqref{eq:UWF}, so in the topology of
spatial distributions the pullback is ordinary reparametrisation.
This identifies $i^*U_{f'}$ with $U_f$ and proves~\eqref{compatcond}. The family $(U_f)$ defines a
tame distribution $U$ on $Q$.

\begin{proposition}\label{prop:WFmetrics}
For $g\in Q$ and $\rho\in\dotT\Sigma_-$ define $L_{g,\rho}\in T_g^*Q$ by
\[
 L_{g,\rho}(h)=\frac12\int_{\gamma_{g,\rho}}
 h(\dot\gamma_{g,\rho},\dot\gamma_{g,\rho})\,\der r,
\]
where $\gamma_{g,\rho}$ is the future-directed ray with initial covector
as in~\eqref{eq:introA}. Then $U$ has a wavefront set in the above
sense, and
\[
 \WF_V(U)=\{(g,\pm L_{g,\rho}):g\in Q,\ \rho\in\dotT\Sigma_-\}.
\]
\end{proposition}

\begin{proof}
Denote the right-hand side by $\mathcal W_0$, and let $h_0$ be the direction
constructed below, so that $L_{g,\rho}(h_0)>0$ for every $g\in Q$ and
every $\rho\in\dotT\Sigma_-$.

Since $A^f_{b,\rho}=(df_b)^*L_{f(b),\rho}$, Theorem~\ref{thm:Sobolev}
gives $\WF_V(U_f)=f^*\mathcal W_0$ for every $f\in\mathcal M$. The fibres of
$\mathcal W_0$ are weak$^*$-closed. To see this, let $c_\alpha L_{g,\rho_\alpha}$, with $c_\alpha>0$ and
$\rho_\alpha$ in the initial unit cosphere bundle, be a net converging
weak$^*$ to some $\zeta\ne0$; nets in the other half of the fibre are
treated in the same way. By compactness we may pass
to a subnet with $\rho_\alpha\to\rho$. Since $\rho\mapsto L_{g,\rho}(h)$
is continuous for every $h$ and $L_{g,\rho}(h_0)>0$, evaluation at $h_0$
shows that $c_\alpha$ converges to some $c\ge0$. Hence
$\zeta=cL_{g,\rho}$, and $c>0$ because $\zeta\ne0$. Thus $\mathcal W_0$ is
admissible.

Now let $\mathcal W$ be any admissible cone, and fix $g\in Q$ and
$\rho\in\dotT\Sigma_-$. For a finite set $\mathcal F$ of directions
containing $h_0$, the affine family
$f_{\mathcal F}(b)=g+\sum_{h\in\mathcal F}b_hh$ takes values in $Q$ for
small $b$ and satisfies Assumption~\ref{ass:A} by its $h_0$ direction.
After shrinking its domain, $f_{\mathcal F}\in\mathcal M$. Its ray
covector at $b=0$ is $(L_{g,\rho}(h))_{h\in\mathcal F}$, which lies in
$\WF_V(U_{f_{\mathcal F}})\subset f_{\mathcal F}^*\mathcal W$. Hence there is
$\xi_{\mathcal F}\in \mathcal W_g$ with $\xi_{\mathcal F}(h)=L_{g,\rho}(h)$ for
every $h\in\mathcal F$. Ordered by inclusion, the net
$(\xi_{\mathcal F})$ converges weak$^*$ to $L_{g,\rho}$, so
$L_{g,\rho}\in \mathcal W_g$. The same argument applied to the ray covector
$-(L_{g,\rho}(h))_{h\in\mathcal F}$ gives $-L_{g,\rho}\in \mathcal W_g$. Thus
$\mathcal W_0\subset \mathcal W$, and $\mathcal W_0$ is the smallest admissible cone.
\end{proof}

If the lower-order terms vary independently, the same argument applies
in the local coordinates $(g,\nabla,W)$ of \eqref{eq:Weitzenboeck},
with the ray covectors acting on the metric component.

To achieve Assumption~\ref{ass:A}, enlarge $K$ to contain a compact
time band between $\Sigma_-$ and $\Sigma_+$ and choose
$h_0=\chi(t)\,\der t\otimes\der t$ supported in this band, with smooth
$\chi\ge0$ positive on an open time interval. Each null ray
crosses this interval and $\der t(\dot\gamma)$ is nonzero there, so
$L_{g,\rho}(h_0)>0$ for every ray. Given any smooth metric family $g_b$,
adjoin the single parameter $a$ by setting $g_{b,a}=g_b+ah_0$.
After shrinking the parameter neighbourhood, this family takes values
in $Q$ and satisfies Assumption~\ref{ass:A} by its $a$ direction.
The original family is its restriction to $a=0$. This supplies the
required local extensions into $\mathcal M$.

\section{Reconstruction of the Cauchy propagation relation}\label{sec:reconstruction}

In this section we show that membership in the Cauchy propagation relation can be tested using diffeomorphisms that displace the two endpoint
base points and pull back the metric. Two signed
normal displacements measure the spatial covector norms; a tangential
displacement adapted to a candidate pair reaches the sum of these norms
only at its prescribed base points and directions. The geometry of the
endpoints and their covectors is illustrated in Figure~\ref{fig:slab}.
This gives the following three-parameter test.

\begin{theorem}[Reconstruction of Cauchy propagation]\label{prop:reconstruction}
Fix a normally hyperbolic operator $P$ on $E$ over $(M,g)$ and two compact
spacelike Cauchy surfaces $\Sigma_\pm$ as above. At $b=0$ write $\mathcal C^\sigma=\mathcal C_0^\sigma$ for its
Cauchy propagation relation on sheet $\sigma$. Fix a candidate pair
$\rho_\pm^0=(z_\pm,\eta_\pm^0)$ and put $a_\pm=|\eta_\pm^0|>0$, using
the induced Riemannian metrics $-g|_{T\Sigma_\pm}$ for spatial norms and
duals. There is a smooth family of compactly supported diffeomorphisms
$F_b$, $b=(b_+,b_0,b_-)\in B\subset\R^3$,
with $F_0=\Id$, such that
\[
 g_b=F_b^*g,\qquad
 P_b=F_b^*P(F_b^*)^{-1}
\]
satisfies Assumption~\ref{ass:A} and, for every $r,s\in\R$,
\begin{equation}\label{eq:reconstruction}
 (\rho_+^0,\rho_-^0)\in\mathcal C^\sigma
 \quad\Longleftrightarrow\quad
 (0;a_+,a_++a_-,a_-)
 \in\WF_{\calL(\calH^r,\calH^s)}(U^\sigma).
\end{equation}
The diffeomorphisms \(F_b\) can be chosen to be the identity outside arbitrarily small neighbourhoods of the two Cauchy surfaces, and their construction requires only the metric near these surfaces and the prescribed endpoint covectors.
Pullback on sections uses a smooth bundle lift. The evolution is taken
between the fixed surfaces $\Sigma_\pm$, and the family is chosen
separately for each candidate pair and sheet. The same test holds with
$U$ in place of $U^\sigma$. For a Dirac-type operator $D$, take
$D_b=F_b^*D(F_b^*)^{-1}$; the assertions hold with
$V^\sigma,V$ and the spaces $\calL(H^r,H^s)$.
\end{theorem}

\begin{proof}
Choose smooth tangential vector fields $W_\pm$ on $\Sigma_\pm$ with
$|W_\pm|\le1$, their norms attaining one only at $z_\pm$. At these points
$W_+$ is the unit Riemannian dual of $\eta_+^0$, while $W_-$ is the
negative unit dual of $\eta_-^0$. Such fields are obtained from local unit
fields and bump functions with unique maxima. In disjoint compact
collars choose $X_\pm$ vanishing near the opposite surface, equal on
their respective surfaces to $-n_+$ and $n_-$ for $\sigma=\mathrm{pos}$,
and to $n_+$ and $-n_-$ for $\sigma=\mathrm{neg}$, where $n_\pm$ are the
future-directed unit normals. Choose $X_0$ restricting to $W_\pm$ at both ends.
Extend the fields with compact support, and let $F_b$ be the time-one
flow of $b_+X_++b_0X_0+b_-X_-$. The flows and bundle lifts depend
smoothly on $b$. Enlarging $K$ and shrinking $B$ gives the support
and common Cauchy temporal hypotheses of Section~\ref{sec:setup}.
Pullback relates the Cauchy data for $P_b$ to data on $F_b(\Sigma_\pm)$
for $P$; changing from the pulled-back connection to the fixed
auxiliary connection adds only a zero-order term.

On the two null sheets (positive and negative, respectively), the normal components of an endpoint covector
are $-|\eta_\pm|$ and $|\eta_\pm|$. The boundary identity
\eqref{eq:boundary} and the parameter-covector formula~\eqref{eq:tau},
with their signed Hamiltonian parameters, consequently give at $b=0$
on the chosen sheet
\[
 (\tau_+,\tau_0,\tau_-)
 =\bigl(|\eta_+|,\ \eta_+(W_+)-\eta_-(W_-),\ |\eta_-|\bigr).
\]
The other sheet has the opposite normal components. In particular
$|\tau_-|=|\eta_-|$ on either sheet, and compactness of the unit
cosphere bundle and smooth dependence give Assumption~\ref{ass:A}
after shrinking $B$.

By Cauchy--Schwarz the middle component satisfies
$\tau_0\le|\eta_+|+|\eta_-|=\tau_++\tau_-$. Equality forces both
endpoints to be $z_\pm$ and their covectors to point in the prescribed
directions $\eta_\pm^0$. Fixing the outer components to be $a_\pm$
also fixes their lengths, so the parameter covector
$(a_+,a_++a_-,a_-)$ occurs exactly at the candidate pair.
The wavefront equality of Theorem~\ref{thm:main} proves
\eqref{eq:reconstruction}. The opposite signs of the normal components
exclude this covector on the other sheet, so the same test works for
unpolarised evolution. The Dirac assertions follow from
Corollary~\ref{cor:introDirac}.
\end{proof}

Note here that the choice of diffeomorphisms and the reconstruction test depend only on the principal symbol of the operator, and not on its lower-order terms.

The ray functional in \eqref{eq:introA} is a symmetric two-tensor light-ray integral, of the kind studied in \cites{FeizmohammadiIlmavirtaOksanenLight,KrishnanSenapatiVashisthLight}. Its conformal and potential-tensor gauge features have a simple meaning here: for $h=fg$ the integral vanishes because $\dot\gamma$ is null, while for $h=\mathcal L_Xg$ it reduces to the endpoint term in \eqref{eq:boundary}. A potential variation with $X$ vanishing at both Cauchy surfaces is therefore invisible to these ray integrals; the endpoint terms survive when $X$ does not vanish there.

\section{Dirac evolution}\label{sec:Dirac}

Let $D_b:C^\infty(M;E)\to C^\infty(M;E)$ be a smooth family of Dirac-type operators on the fixed bundle $E$. We assume that $D_b=D_0$ outside $K$ and that $P_b=D_b^2$ has principal symbol $-g_b^{-1}\Id_E$. Then $P_b$ is a family of normally hyperbolic operators satisfying the assumptions of Section~\ref{sec:setup}. By the Schr\"odinger--Lichnerowicz formula for the square of a geometric Dirac operator, the Weitzenb\"ock connection in~\eqref{eq:Weitzenboeck} is the spin connection; in the twisted case it is the spin connection tensored with the twisting connection.

\begin{remark}
Geometric Dirac operators associated with a smooth family of metrics on a spin manifold are included by identifying the spinor bundles locally in parameter space. The standard identification, chosen to be the identity outside $K$, pulls the operators back to a fixed bundle; see \cite{BGM}*{Section 5}. The same applies to twisted spinor bundles.
\end{remark}

\subsection{The two Cauchy problems}

Let $D_b$ be a family as in Corollary~\ref{cor:introDirac}. In a collar of a spacelike Cauchy surface $\Sigma$ we write
\begin{equation}\label{eq:Diracnormal}
 D_b=a_{b,\Sigma}\nabla_{n_b}+T_{b,\Sigma},
 \qquad L_{b,\Sigma}=-a_{b,\Sigma}^{-1}T_{b,\Sigma},
\end{equation}
where $T_{b,\Sigma}$ is tangential of order one, including its zero-order term. The coefficient $a_{b,\Sigma}$ is invertible since the normal covector is noncharacteristic. The Dirac equation implies
\[
 \nabla_{n_b}u|_\Sigma=L_{b,\Sigma}(u|_\Sigma).
\]
Thus the second component of the Cauchy data is determined by the first. Define
\begin{equation}\label{eq:DiracJ}
 J_{b,\Sigma}f=\begin{pmatrix}f\\L_{b,\Sigma}f\end{pmatrix}.
\end{equation}

\begin{proposition}\label{prop:Diracdata}
The solution of the second-order Cauchy problem
\[
 D_b^2u=0,\qquad R_{b,\Sigma}u=J_{b,\Sigma}f
\]
is the unique solution of $D_bu=0$ with $u|_\Sigma=f$.
\end{proposition}

\begin{proof}
A solution of $D_bu=0$ solves $D_b^2u=0$ and has the stated Cauchy data. Conversely, let $u$ solve the second-order problem and set $w=D_bu$. Then $P_bw=D_bP_bu=0$ and, by \eqref{eq:Diracnormal}, $w|_\Sigma=0$. Moreover, $D_bw=P_bu=0$, so the same formula gives $\nabla_{n_b}w|_\Sigma=0$. Uniqueness for the second-order Cauchy problem implies $w=0$. Thus existence and uniqueness for the Dirac problem follow from those for its square.
\end{proof}

Write $L_{b,\pm}=L_{b,\Sigma_\pm}$ and $J_{b,\pm}=J_{b,\Sigma_\pm}$. These operators depend smoothly on $b$. Let $\operatorname{pr}_1$ be projection onto the first component. If $U_b$ is the Cauchy evolution for $D_b^2$, the Dirac evolution is therefore
\begin{equation}\label{eq:DiracU}
 V_b=\operatorname{pr}_1U_bJ_{b,-},\qquad U_bJ_{b,-}=J_{b,+}V_b.
\end{equation}
The second identity expresses invariance of the constrained Cauchy data. For $s\in\R$, the maps
\[
 J_{b,\pm}:H^s(\Sigma_\pm;E_\pm)\longrightarrow\calH^{s-\frac12}(\Sigma_\pm),
 \qquad
 \operatorname{pr}_1:\calH^{s-\frac12}(\Sigma_\pm)\longrightarrow H^s(\Sigma_\pm;E_\pm)
\]
are bounded. Thus \eqref{eq:DiracU} also gives boundedness of $V_b$ on each fixed Sobolev scale. Reversing the Cauchy surfaces gives its bounded inverse.

\subsection{Green operators and the wavefrontset inclusions}

Let $G_{b,\mathrm{ret/adv}}$ be the Green operators for $P_b=D_b^2$. The Green operators for $D_b$ are
\begin{equation}\label{eq:DiracGreen}
 S_{b,\mathrm{ret/adv}}
 =D_bG_{b,\mathrm{ret/adv}}
 =G_{b,\mathrm{ret/adv}}D_b.
\end{equation}
For $f\in C^\infty_0(M;E)$ the difference
$D_bG_{b,\mathrm{ret/adv}}f-G_{b,\mathrm{ret/adv}}D_bf$ solves the homogeneous equation for $P_b$ and has the corresponding causal support. Uniqueness implies that it vanishes. Consequently
\[
 D_bS_{b,\mathrm{ret/adv}}f=f,
 \qquad S_{b,\mathrm{ret/adv}}D_bf=f,
\]
and differentiation does not enlarge support. Taking the difference of the retarded and advanced Green operators gives
\begin{equation}\label{eq:Diraccausal}
 S_b=S_{b,\mathrm{ret}}-S_{b,\mathrm{adv}}=D_bG_b.
\end{equation}

There is also a useful formula in terms of the second-order Cauchy solution operator. Denote by $\mathcal T_{b,\Sigma}$ the solution operator for
\[
 P_bw=0,\qquad w|_\Sigma=0,\qquad\nabla_{n_b}w|_\Sigma=f.
\]
Then $D_b\mathcal T_{b,\Sigma}a_{b,\Sigma}^{-1}f$ solves $D_bu=0$ and restricts to $f$ at $\Sigma$. Consequently,
\begin{equation}\label{eq:DiracCauchyGreen}
 V_b=\operatorname{res}_{\Sigma_+}D_b\mathcal T_{b,\Sigma_-}a_{b,-}^{-1},
 \qquad a_{b,-}=a_{b,\Sigma_-}.
\end{equation}
This is the same construction used in \cite{BaerStrohmaierIndex}*{Section 2} to characterise the Cauchy evolution for the Dirac equation as a Fourier integral operator.
Together with \eqref{eq:DiracGreen}, it expresses the first-order solution operators in terms of those of the square.

Since $D_b$ is a differential operator with coefficients smooth in $(b,x)$, applying it in the output variable does not enlarge the joint wavefront set. Now Theorem~\ref{thm:kernel} gives
\[
 \WF(K_{S_{\mathrm{ret/adv}}})\subset\Lambda'\cup\Delta^*,
 \qquad \WF(K_S)\subset\Lambda'.
\]
Thus passage from the square to the first-order operator introduces no additional parameter covectors.

\begin{proof}[Proof of Corollary~\ref{cor:introDirac}]
For fixed $r,s$, the continuous maps
\[
 A\longmapsto\operatorname{pr}_1AJ_{b,-}
\]
from $\calL(\calH^{r-\frac12}(\Sigma_-),\calH^{s-\frac12}(\Sigma_+))$ to
$\calL(H^r(\Sigma_-;E_-),H^s(\Sigma_+;E_+))$ depend smoothly on $b$ in operator norm. Applying such a smooth operator family to a Banach-valued distribution does not enlarge its wavefront set. Theorem~\ref{thm:main} and~\eqref{eq:DiracU} therefore give the bound for $V$.
For the two sheets put
\begin{equation}\label{eq:Diracpolarized}
 V_b^\sigma=\operatorname{pr}_1U_b\Pi_b^\sigma J_{b,-}.
\end{equation}
They sum to $V_b$ modulo a jointly smooth kernel and satisfy the polarised inclusions by the same argument. To prove the reverse inclusions, let $l_b(\eta)$ be the principal symbol of $L_{b,\Sigma}$ in~\eqref{eq:Diracnormal}. At a point of $\Sigma$, in coordinates adapted to its unit normal, the Clifford relations give
\[
 l_b(\eta)^2=-|\eta|_b^2\Id,\qquad
 a_{b,\Sigma}l_b(\eta)=-l_b(\eta)a_{b,\Sigma},
\]
where $|\eta|_b$ is the norm for the induced Riemannian metric. For $\eta\ne0$, the eigenspaces with eigenvalues $\pm\rmi|\eta|_b$ have equal positive rank, since the invertible coefficient $a_{b,\Sigma}$ interchanges them. Write $e_b^\sigma$ for the projection onto the eigenspace corresponding to sheet $\sigma$. The principal symbol of $\Pi_b^\sigma J_{b,-}$ is the Cauchy graph inclusion composed with $e_b^\sigma$.

Transport for $D_b^2$ is an isomorphism between the initial and final
Cauchy polarisations on each sheet. Projection of a Cauchy polarisation
onto its first component is also an isomorphism, since the polarisation
is a graph. The principal symbol of $V^\sigma$ is therefore the nonzero
projection $e_b^\sigma$ followed by these isomorphisms, and does not vanish
on $(\Lambda_U^\sigma)'$.

To exclude cancellation in the sum, let $\mathcal U_b$ be the spacetime
Cauchy solution operator for $D_b^2$ with initial surface $\Sigma_-$,
and put $E_b^\sigma=\operatorname{pr}_1\Pi_b^\sigma J_{b,-}$. The
$E_b^\sigma$ are smooth families of pseudodifferential operators of order
zero. The kernels of
\[
 D_b\mathcal U_b\Pi_b^\sigma J_{b,-}
\]
have wavefront sets on disjoint spacetime null sheets. Their sum is
jointly smooth, since $D_b\mathcal U_bJ_{b,-}=0$ and the projections sum
to the identity modulo a jointly smooth kernel. Hence each residual is
jointly smooth. Restricting to $\Sigma_-$ and using
\eqref{eq:Diracnormal} gives
\[
 \Pi_b^\sigma J_{b,-}=J_{b,-}E_b^\sigma
\]
modulo a jointly smooth kernel. We then have $V_b^\sigma=V_bE_b^\sigma$
modulo such a kernel. Right composition is microlocal at the nonzero input covectors in these relations, so
\[
 \WF(K_{V^\sigma})=(\Lambda_U^\sigma)',\qquad
 \WF(K_V)=(\Lambda_U^{\mathrm{pos}})'\cup(\Lambda_U^{\mathrm{neg}})'.
\]
Lemma~\ref{lem:kernel} now gives the asserted operator-valued equalities.
At kernel level, composition with $J_{b,-}$ consists of tangential differential operators with coefficients smooth in $b$. It preserves the smooth-kernel-valued distribution property, which proves the averaging statement. The polar-cone inclusions follow from $\Gamma_{\mathrm{pos}}\subset C^\circ$ and $\Gamma_{\mathrm{neg}}\subset-C^\circ$.
\end{proof}

\appendix

\section{Construction of the parameter-dependent fundamental solutions}\label{app:kernel}

Here we give a detailed proof of Theorem~\ref{thm:kernel} essentially following the construction of \cite{StrohmaierStress}*{Appendix C} for the more general case of vector bundles. This is a parameter-dependent version of the parametrix construction of \cites{MR0388464,MR1362544}.

\subsection{Construction and normalisation of the parametrix}

\begin{proposition}\label{prop:jointparametrix}
There is a bundle-valued Lagrangian distribution
\[
 \widehat G\in
 I^{-\frac32-\frac{\dim B}{4}}
 (B_0\times(I\times\Sigma)^2,\Lambda';E\boxtimes E^*)
\]
such that $F_b=P_b\widehat G_b$ is jointly smooth and vanishes to infinite
order on the equal-time hypersurface
\[
 \mathcal S=\{(b,t,x,s,y):t=s\}.
\]
Moreover,
\begin{equation}\label{eq:exactequaldata}
 \widehat G_b(s,x;s,y)=0,\qquad
 \left.\partial_t\widehat G_b(t,x;s,y)\right|_{t=s}
       =c_b(s,x)^{-1}\delta_\Sigma(x,y).
\end{equation}
\end{proposition}

\begin{proof}
We first construct a Lagrangian parametrix $\widetilde G$ with the
identities~\eqref{eq:exactequaldata} holding modulo jointly smooth kernels.
We then correct its full equal-time jet by a smooth kernel.

Put $m=-\tfrac32-\tfrac{\dim B}{4}$ and denote the two null sheets of
$\Lambda'$ by $\Lambda'_\pm$. The scalar principal symbol of $P_b$ in
the output variable is $p=-2H_b$. It vanishes on $\Lambda'_\pm$, since
the output covectors are null, and $\der p$ is nonzero there. Applying
the product formula on $B_0\times(I\times\Sigma)^2$, with the parameters
included among the base variables, gives, for
$A\in I^q(\Lambda'_\pm;E\boxtimes E^*)$,
\begin{equation}\label{eq:bundletransport}
 P_bA\in I^{q+1}(\Lambda'_\pm;E\boxtimes E^*),\qquad
 \sigma_{q+1}(P_bA)=\mathcal T\sigma_q(A).
\end{equation}
For this formula and the successive correction below see
\cite{MR0388464}*{Theorems 5.3.1 and 5.3.2} and
\cite{MR2512677}*{Theorem 25.2.4 and Lemma 26.1.16}; the bundle-valued
versions are given in \cite{IslamStrohmaierFeynman}*{Theorem E.1 and Appendix F}.

In a local bundle frame and half-density notation,
$\mathcal T=-\rmi\mathcal L_{X_p}+p_{\mathrm{sub}}$, where
$p_{\mathrm{sub}}$ is the subprincipal symbol. The combination is
independent of these choices. By~\eqref{eq:Weitzenboeck}, its bundle part
is the covariant derivative for the Weitzenb\"ock connection of $P_b$;
see \cite{IslamStrohmaierFeynman}*{Theorem 2.8}. The homogeneous transport
equation is therefore parallel transport for this connection, together
with the natural transport of the half-density and Maslov factors.
On the joint cotangent bundle $X_p=-2X_H$, including the parameter
covector component. Its integral curves on $\Lambda'_\pm$ are thus the
lifted null bicharacteristics in~\eqref{eq:tau}, with reversed orientation
and rescaled parameter.

Since $t=s$ is a Cauchy surface, each integral curve meets it exactly
once with the input variables fixed. The intersection is transverse
because the time component of $X_p$ is nonzero. Solving $\mathcal Ta=f$
with prescribed initial symbol at $t=s$ is therefore a linear ordinary
differential equation along each ray, in the homomorphism bundle tensored
with the Maslov and half-density lines. Its solution is unique and smooth
in the initial variables and parameters.

At equal times, pullback of either sheet gives the spatial diagonal,
with zero covectors in the $b$ and $s$ variables. The trace symbol
identification, using the fixed density $\der\mu$, identifies initial
symbols with endomorphism-valued symbols on $\dotT\Sigma$. For our choice
of $m$, the equal-time trace of an element of $I^{m-j}(\Lambda'_\pm)$
has tangential pseudodifferential order $-1-j$, and its $t$ derivative
has order $-j$. If $\alpha^\pm$ are the initial trace symbols, the
principal symbols of the two traces of their sum are
\[
 \alpha^++\alpha^-,\qquad
 \rmi\lambda_b^+\alpha^++\rmi\lambda_b^-\alpha^-,
\]
where $\lambda_b^+>\lambda_b^-$ are the two roots of
$g_b^{-1}(\lambda\der t+\eta,\lambda\der t+\eta)=0$.
They are smooth and homogeneous of degree one for $\eta\ne0$, and their
difference is bounded below by a positive multiple of $|\eta|$ on
compact sets.

Choose the leading symbols by solving $\mathcal Ta_0^\pm=0$ with initial
values
\[
 \alpha_0^+=
 \frac{c_b(s)^{-1}}{\rmi(\lambda_b^+-\lambda_b^-)}\Id,
 \qquad \alpha_0^-=-\alpha_0^+.
\]
Surjectivity of the principal symbol map gives representatives
$A_0^\pm\in I^m(\Lambda'_\pm;E\boxtimes E^*)$.
Their sum satisfies the equation and both Cauchy identities to leading
order. We now remove the lower-order errors successively.

Suppose $A^{(j-1)}=\sum_{k<j}(A_k^++A_k^-)$ has been constructed, with
$A_k^\pm\in I^{m-k}(\Lambda'_\pm)$, so that its equation error has order
$m+1-j$ and its two Cauchy errors have tangential orders $-1-j$ and $-j$.
Let $r_j^\pm$ be the principal symbols of the equation error on the two
sheets, and let $e_{0,j},e_{1,j}$ be those of
$A^{(j-1)}|_{t=s}$ and
$\partial_tA^{(j-1)}|_{t=s}-c_b(s)^{-1}\delta_\Sigma$.
Solve
\begin{equation}\label{eq:transportrecursion}
 \mathcal Ta_j^\pm=-r_j^\pm
\end{equation}
with initial values determined by
\begin{equation}\label{eq:amplitudeinitial}
 \begin{pmatrix}1&1\\ \rmi\lambda_b^+&\rmi\lambda_b^-\end{pmatrix}
 \begin{pmatrix}\alpha_j^+\\\alpha_j^-\end{pmatrix}
 =-\begin{pmatrix}e_{0,j}\\e_{1,j}\end{pmatrix}.
\end{equation}
The determinant is $\rmi(\lambda_b^--\lambda_b^+)$, so the initial symbols
exist and have degree $-1-j$. Choose representatives
$A_j^\pm\in I^{m-j}(\Lambda'_\pm)$ with these principal symbols.
By~\eqref{eq:bundletransport} and the trace calculus, adding them lowers
all three errors by one order. 

Asymptotic summation gives
$\widetilde G\sim\sum_{j\ge0}(A_j^++A_j^-)$ and hence
\begin{equation}\label{eq:parametrixremainders}
 P_b\widetilde G_b\in C^\infty,
 \qquad \widetilde G_b|_{t=s}\in C^\infty,
 \qquad \partial_t\widetilde G_b|_{t=s}
                    -c_b(s)^{-1}\delta_\Sigma\in C^\infty.
\end{equation}
Differentiating the transport equations gives symbol
estimates for all derivatives in the parameters and initial variables,
uniformly on compact parameter sets. Asymptotic summation with these estimates gives the
jointly smooth remainders above.

It remains to impose the equal-time identities to infinite order. This
step can be carried out using formal jets computed entirely from the
coefficients of $P_b$. Define distribution kernels $J_{b,j}(s)$ by
\[
 J_{b,0}=0,\qquad J_{b,1}=c_b(s)^{-1}\delta_\Sigma,
\]
and recursively, for $k\ge0$, by
\begin{align}\label{eq:formaljets}
 c_bJ_{b,k+2}
 ={}&-\sum_{\ell=1}^k\binom{k}{\ell}
               (\partial_t^\ell c_b)J_{b,k-\ell+2}\notag\\
    &-\sum_{\ell=0}^k\binom{k}{\ell}
               (\partial_t^\ell A_{b,1})J_{b,k-\ell+1}
     -\sum_{\ell=0}^k\binom{k}{\ell}
               (\partial_t^\ell A_{b,2})J_{b,k-\ell}.
\end{align}
All coefficients on the right are evaluated at $t=s$ and act in the
output variable. These are the identities obtained by differentiating
$P_bG_b=0$ formally at equal times. In particular, each $J_{b,j}$ is the
kernel of a differential operator with coefficients smooth in $(b,s)$.

The restrictions $\partial_t^j\widetilde G_b|_{t=s}$ exist. Indeed, a
nonzero covector in the conormal bundle of $\mathcal S$ is timelike in
each spacetime factor, whereas the corresponding covectors of $\Lambda'$
are null. Put
\[
 e_{b,j}=\left.\partial_t^j\widetilde G_b\right|_{t=s}-J_{b,j}.
\]
The first two kernels $e_{b,0},e_{b,1}$ are jointly smooth
by~\eqref{eq:parametrixremainders}. Applying $\partial_t^k$ to
$P_b\widetilde G_b$ and restricting to $t=s$ gives the
recursion~\eqref{eq:formaljets} for these differences, with the additional
jointly smooth term $\partial_t^k(P_b\widetilde G_b)|_{t=s}$. Since $c_b$
is invertible, induction proves that every $e_{b,j}$ is jointly smooth.

Apply Borel's extension theorem to these smooth bundle-valued jets on
$\mathcal S$, with $(b,s,x,y)$ as tangential variables and $q=t-s$ as the
normal variable. To see the parameter dependence explicitly, in a
trivialisation one may take a sum of the terms
\[
 \chi_0(q/\varepsilon_j)\frac{q^j}{j!}e_{b,j}(s,x,y),
\]
where $\chi_0=1$ near zero and the positive numbers $\varepsilon_j$
decrease sufficiently fast to control successively all derivatives on a
compact exhaustion. A partition of unity on $\mathcal S$ gives a smooth
global kernel $\mathcal B$ with $\partial_t^j\mathcal B|_{t=s}=e_{b,j}$ for every $j$; cf.
\cite{Ho1}*{Theorem 1.2.6}. 

Set $\widehat G=\widetilde G-\mathcal B$. It has all the prescribed jets
$J_{b,j}$ and remains in the same Lagrangian class. Its residual
$F_b=P_b\widehat G_b$ is jointly smooth. The
recursion~\eqref{eq:formaljets} implies
$\partial_t^kF_b|_{t=s}=0$ for every $k$. Since $\partial_t$ is transverse
to $\mathcal S$, this is exactly infinite order vanishing along
$\mathcal S$. The proposition follows.
\end{proof}

Define the retarded and advanced parametrices by
\begin{equation}\label{eq:Heavisideparametrices}
 T_{b,\mathrm{ret}}=\theta(t-s)\widehat G_b,
 \qquad T_{b,\mathrm{adv}}=-\theta(s-t)\widehat G_b.
\end{equation}
These products are well-defined by the wavefront criterion for products, \cite{Ho1}*{Theorem 8.2.10}:
the conormal covectors of $t=s$ cannot cancel the nonzero null covectors
of $\Lambda'$. The equal-time identities give
\begin{equation}\label{eq:retardederror}
 P_bT_{b,\mathrm{ret}}=\Id+R_{b,\mathrm{ret}},\qquad
 P_bT_{b,\mathrm{adv}}=\Id+R_{b,\mathrm{adv}},
\end{equation}
where
\[
 R_{b,\mathrm{ret}}=\theta(t-s)F_b,
 \qquad R_{b,\mathrm{adv}}=-\theta(s-t)F_b.
\]
For example, since $\widehat G_b|_{t=s}=0$, the term
$c_b(t)\delta'(t-s)\widehat G_b$ equals
$-c_b(s)\delta(t-s)\partial_t\widehat G_b|_{t=s}$. This cancels one
of the two cross terms $c_b(s)\delta(t-s)\partial_t\widehat G_b|_{t=s}$
from differentiating $\theta(t-s)\widehat G_b$ twice. The remaining
boundary term is the identity kernel.
Both remainders are jointly smooth because $F$
is flat on $\mathcal S$. Their supports satisfy $t\ge s$ and $t\le s$,
respectively.

We also record the joint wavefront set inclusion
\begin{equation}\label{eq:HeavisideWF}
 \WF(T_{\mathrm{ret/adv}})\subset\Lambda'\cup\Delta^*.
\end{equation}
Away from $\mathcal S$ it follows directly from the bound for
$\widehat G$. On $\mathcal S$ away from $x=y$, that kernel is smooth and
all its normal jets vanish, since each $J_{b,j}$ is supported on the
spatial diagonal. Thus multiplication by the Heaviside function is smooth
there. On the spatial diagonal the covectors of $\Lambda'$ have parameter
component zero and opposite spacetime components. Adding a conormal
covector of $t=s$ retains this property. All additional covectors given by
the product estimate consequently belong to $\Delta^*$. This proves
\eqref{eq:HeavisideWF}.

\subsection{The Neumann series}

We give the convergence argument for the retarded remainder. The same
argument with time reversed treats the advanced remainder.

\begin{lemma}\label{lem:Volterra}
Let $R_b(t,x;s,y)$ be a jointly smooth homomorphism-valued kernel on a
neighbourhood of $Q\times([a,c]\times\Sigma)^2$, where $Q\Subset B$ is
compact, and suppose it vanishes for $t<s$. Then
\begin{equation}\label{eq:Neumannkernel}
 Q_b^{\mathrm{ret}}=\sum_{j=1}^\infty(-1)^jR_b^j
\end{equation}
converges in the $L^2$ operator norm and with all spacetime derivatives,
uniformly for $b\in Q$. The sum is jointly smooth in the spacetime and
parameter variables, has support in $t\ge s$, and satisfies
\[
 (\Id+R_b)^{-1}=\Id+Q_b^{\mathrm{ret}}.
\]
The inverse identity holds on $L^2$ and on continuous sections on the
slab, and also on smooth sections depending smoothly on additional
parameters.
\end{lemma}

\begin{proof}
Put $Z=[a,c]\times\Sigma$, fix an auxiliary Hermitian bundle metric,
and use the product density fixed above. Let $V=\mu(\Sigma)$, $L=c-a$,
and let $M$ be a uniform bound on the kernel of $R_b$. The ordered time
simplex in the kernel of $R_b^j$ has volume
$(t-s)^{j-1}/(j-1)!$. Consequently, for $t\ge s$, its kernel is bounded
by $M^jV^{j-1}(t-s)^{j-1}/(j-1)!$. Squaring and integrating over
$Z\times Z$ gives the Hilbert--Schmidt estimate
\begin{equation}\label{eq:Volterraestimate}
 \begin{aligned}
 \|R_b^j\|_{\calL(L^2(Z;E))}
 &\le\|R_b^j\|_{\mathrm{HS}}\le\frac{(MVL)^j}{(j-1)!\sqrt{2j(2j-1)}}\le\frac{(MVL)^j}{j!}.
 \end{aligned}
\end{equation}
Thus $S_b=\sum_{j\ge0}(-R_b)^j$ converges in operator norm, uniformly on
$Q$. Multiplying its partial sums by $\Id+R_b$ shows that
$S_b=(\Id+R_b)^{-1}$.

For every bounded operator $T$ on $L^2(Z;E)$, the two smoothing factors
give a smooth kernel and the estimates
\[
 \|R_bTR_b\|_{C^k(Z\times Z)}\le C_k\|T\|_{\calL(L^2(Z;E))},
\]
uniformly on $Q$. Indeed, derivatives in the two kernel variables act
only on the two outer factors, and the estimate follows from
the Cauchy--Schwarz inequality. For $N\ge2$, the partial sum of~\eqref{eq:Neumannkernel}
is
\[
 -R_b+R_b\left(\sum_{j=0}^{N-2}(-R_b)^j\right)R_b.
\]
It therefore converges with all spacetime derivatives to
$Q_b^{\mathrm{ret}}=-R_b+R_bS_bR_b$. The support assertion follows
from the support of every term.

Finally, $b\mapsto R_b$ is smooth as a function taking values in the Banach space of bounded operators in $L^2$. 
Smoothness of
inversion gives smooth dependence of $S_b$; in particular,
$\partial_{b^\nu}S_b=-S_b(\partial_{b^\nu}R_b)S_b$.
The same smoothing-factor estimates applied to $-R_b+R_bS_bR_b$ and its
parameter derivatives show that its kernel is jointly smooth. The
inverse identities consequently hold on continuous and smooth sections,
also with smooth auxiliary parameters.
\end{proof}

\subsection{Correction to the fundamental solutions}

\begin{proof}[Proof of Theorem~\ref{thm:kernel}]
Use Proposition~\ref{prop:jointparametrix} and
\eqref{eq:Heavisideparametrices} on a time interval slightly larger than
the compact slab under consideration. Apply Lemma~\ref{lem:Volterra} to
$R_{b,\mathrm{ret}}$ and denote the resulting smooth kernel by
$Q_{b,\mathrm{ret}}$. Put
\begin{equation}\label{eq:correctedretarded}
 \mathcal G_{b,\mathrm{ret}}
       =T_{b,\mathrm{ret}}(\Id+Q_{b,\mathrm{ret}}).
\end{equation}
All compositions are well-defined on compact subslabs. The time support
ensures that an intermediate time lies between the input and output
times, and $\Sigma$ is compact. A cutoff in the integration variable may
therefore be inserted when discussing compact subsets of the kernel.
Equation~\eqref{eq:retardederror} and the inverse identity give
\[
 P_b\mathcal G_{b,\mathrm{ret}}=\Id.
\]

The correction $T_{\mathrm{ret}}Q_{\mathrm{ret}}$ has a jointly smooth
kernel. Indeed, neither $\Lambda'$ nor $\Delta^*$ contains a covector
with zero input spacetime component. By the kernel wavefront calculus,
\eqref{eq:HeavisideWF} therefore implies that $T_{\mathrm{ret}}$ maps
jointly smooth forcing with compact support to jointly smooth sections.
The same statement holds with additional smooth variables present.
Apply it to the smooth kernel of $Q_{\mathrm{ret}}$ in its intermediate
variable, using the cutoff just described. This proves
\begin{equation}\label{eq:smoothGreenCorrection}
 \mathcal G_{\mathrm{ret}}-T_{\mathrm{ret}}\in C^\infty,
 \qquad
 \WF(\mathcal G_{\mathrm{ret}})\subset\Lambda'\cup\Delta^*.
\end{equation}

For each fixed $b$ and compactly supported smooth forcing $f$, the section
$\mathcal G_{b,\mathrm{ret}}f$ solves $P_bu=f$ and vanishes before the
support of $f$. Its difference from the usual retarded solution solves
the homogeneous equation with zero Cauchy data to the past of that
support. Uniqueness for the ordinary Cauchy problem identifies
$\mathcal G_{b,\mathrm{ret}}$ with $G_{b,\mathrm{ret}}$ on the slab.
Existence and uniqueness at a fixed parameter are the standard results
for normally hyperbolic bundle operators, see \cite{BGP}*{Chapter 3}.
In particular the corrected operator has the full retarded causal support.

Repeating the construction with time reversed gives
\[
 G_{\mathrm{adv}}-T_{\mathrm{adv}}\in C^\infty,
 \qquad \WF(G_{\mathrm{adv}})\subset\Lambda'\cup\Delta^*.
\]
Since $T_{\mathrm{ret}}-T_{\mathrm{adv}}=\widehat G$, subtraction yields
\[
 G_{\mathrm{ret}}-G_{\mathrm{adv}}-\widehat G\in C^\infty.
\]
Proposition~\ref{prop:jointparametrix} now gives the claimed classical
Lagrangian order of the exact causal propagator. Compact subsets of
$B\times M\times M$ are covered by finitely many parameter neighbourhoods
and finite time slabs of the preceding type. Testing any of the constructed
kernels against smooth compactly supported sections in both spacetime
variables gives a smooth function of $b$, since their wavefront sets have
no nonzero covectors purely in the parameter directions. Thus uniqueness
of the fixed-parameter Green operators on overlaps also identifies the
joint distributions there. This proves the theorem globally.
\end{proof}

\section*{Use of AI tools}

The mathematical statements and proofs in this paper were developed by the author, with the exception of the example in Section~\ref{sec:counterexample}, which emerged in a discussion between the author and GPT-6.
AI tools (GPT-6, Claude Opus 5.5) were used for typesetting, writing assistance, proofreading, and for generating Fig.~\ref{fig:slab}. 
AI was not integral to the research methodology in this work. The author takes full responsibility for the mathematical contents of the paper.

\enlargethispage{6\baselineskip}

\enlargethispage{2\baselineskip}
\end{document}